\documentclass[reqno,11pt]{amsart} 
\usepackage{amsmath} 
\usepackage{amssymb}
\def\b{\beta}

\def\d{\delta} 
 
\def\l{\lambda} 

\def\t{\tau}

\def\L{\Lambda}

\newfam\Bbbfam 
\font\tenBbb=msbm10 
\font\sevenBbb=msbm7 
\font\fiveBbb=msbm5 
\textfont\Bbbfam=\tenBbb 
\scriptfont\Bbbfam=\sevenBbb 
\scriptscriptfont\Bbbfam=\fiveBbb

\newcommand{\R}     {\mathbb{R}} 
\newcommand{\Z}     {\mathbb{Z}} 
\newcommand{\N}     {\mathbb{N}} 
\renewcommand{\P}   {\mathbb{P}} 
 
\newcommand{\E}     {\mathbb{E}} 
\newcommand{\Q}     {\mathbb{Q}}

\def\1{{\mathchoice {1\mskip-4mu\mathrm l}      
{1\mskip-4mu\mathrm l} 
{1\mskip-4.5mu\mathrm l} {1\mskip-5mu\mathrm l}}} 
\newcommand{\ssup}[1] {{\scriptscriptstyle{({#1}})}} 
\def\comment#1{} 
\newtheoremstyle{thm}{2ex}{2ex}{\itshape\rmfamily}{} 
{\bfseries\rmfamily}{}{1.7ex}{} 
 
\newtheoremstyle{rem}{1.3ex}{1.3ex}{\rmfamily}{} 
{\itshape\rmfamily}{}{1.5ex}{} 
 
\newenvironment{proofsect}[1] 
{\vskip0.1cm\noindent{\bf #1.}\hskip0.5cm}

\newtheorem{theorem}{Theorem}[section] 
\newtheorem{lemma}[theorem]{Lemma} 
\newtheorem{prop}[theorem] {Proposition} 
\newtheorem{cor}[theorem]  {Corollary} 
\newtheorem{remark}[theorem]  {Remark} 
\newtheorem{defn}[theorem] {Definition}

\theoremstyle{definition}
\newtheorem{example}[theorem] {Example}

 \newcommand{\s}{\sigma}

\renewcommand{\section}{\secdef\sct\sect} 
\newcommand{\sct}[2][default]{\refstepcounter{section} 
\vspace{0.8cm} 
\setcounter{equation}{0} 
\centerline{ 
\large\scshape \arabic{section}.\ #1} 
\vspace{0.2cm}} 
\newcommand{\sect}[1]{ 
\vspace{0.8cm} 
\centerline{\large\scshape #1} 
\vspace{0.2cm}} 
 
\renewcommand{\subsection}{\secdef \subsct\sbsect} 
\newcommand{\subsct}[2][default]{\refstepcounter{subsection} 
\nopagebreak 
\vspace{0.5\baselineskip} 
{\flushleft\bf \arabic{section}.\arabic{subsection}~\bf #1  } 
\nopagebreak} 
\newcommand{\sbsect}[1]{\vspace{0.1cm}\noindent 
{\bf #1}\vspace{0.1cm}}

\renewcommand{\subsubsection}{%
\secdef \subsubsect\sbsbsect} 
\newcommand{\subsubsect}[2][default]{%
\refstepcounter{subsubsection} 
\nopagebreak 
\vspace{0.1\baselineskip} 
\nopagebreak 
{\flushleft 
\sffamily\slshape 
\arabic{section}.\arabic{subsection}.\arabic{subsubsection} 
\ %
\sffamily #1\/.}\ } 
\newcommand{\sbsbsect}[1]{\vspace{0.1cm}\noindent 
{\bf #1}\ } 
 
\renewcommand{\d}{{\rm d}}

\newcommand{\Leb}{{\rm Leb}} 
\newcommand{\Sym}{\mathfrak{S}}

\newcommand{\dist}{{\operatorname {dist}}}

\newcommand{\Bcal}   {{\mathcal B }}

\newcommand{\Ccal}   {{\mathcal C }} 
 
\newcommand{\Fcal}   {{\mathcal F }} 
 
\newcommand{\Hcal}   {{\mathcal H }} 
\newcommand{\Ical}   {{\mathcal I }}

\newcommand{\Mcal}   {{\mathcal M }}

\newcommand{\Scal}   {{\mathcal S }}

\newcommand{\Xcal}   {{\mathcal X }} 
\newcommand{\Ycal}   {{\mathcal Y }} 
\newcommand{\Zcal}   {{\mathcal Z }} 
\newcommand{\m} {{\mathfrak m}}
 
\begin{document}
\title{Large systems of symmetrized trapped Brownian Bridges and Schr\"odinger processes} 

\maketitle

\centerline {\sc By Stefan Adams ({\textdagger})\/ and Spyros Garouniatis \footnote{ Mathematics Institute, University of Warwick, Coventry CV4 7AL, United Kingdom, {\tt spyridon.garouniatis@warwick.ac.uk}}}
\bigskip
\centerline{November 27 2024}
\bigskip
\begin{abstract}
    Consider a large system of $N$ Brownian motions in $\R ^d$ fixed on a time interval $[0,\beta]$ with symmetrized initial and terminal conditions, under the influence of a trap potential. Such systems describe systems of bosons at positive temperatures confined in a spatial domain. We describe the large $N$ behavior of the averaged path (that is, their empirical path measure) and its connection with a well known optimal transport problem formulated by Erwin Schr\"odinger. We also explore the asymptotic behavior of the Brownian motions in terms of Large Deviations. In particular, the rate function that governs the mean of occupation measures turns out to be the well-known Donsker-Varadhan rate function. We therefore prove a simple formula for the large $N$ asymptotic of the symmetrized trace of $e^{-\beta \Hcal_N}$, where $\Hcal_N$ is an $N$ particle Hamilton operator in a trap
\end{abstract}\bigskip

\noindent \textbf{Acknowledgments}. This article is dedicated to the memory of Stefan Adams who sadly passed away in May 2024. His contribution to this article was invaluable, and his insights and dedication will be greatly missed. \bigskip

\section{Introduction}

\noindent In this paper we study a model that examines how a large system of particles at positive temperature behaves when the particles are trapped and influenced by specific constraints. When this system gets very large, the average path behavior of the particles becomes such that we can approximate it by a single diffusion. The main objective of this paper is to contribute to a probabilistic interpretation of the Bose-Einstein condensation. In Chapter 4 we will further elaborate on the relevant literature and future directions toward this goal. Let us give a brief overview of the model and its connection to quantum statistical mechanics. 

\noindent We consider a system of $N$ non-interacting bosons in a trap potential $ W $. The system is described by the Hamilton operator, defined on symmetric functions of $L^2(\R ^ {Nd})$
\begin{equation}\label{Hdef}
\Hcal_N=\sum_{i=1}^N \big(-\Delta_i+W(x_i)\big),\qquad x_1,\dots,x_N\in\R^d,
\end{equation}
where the $i$-th Laplace operator, $ \Delta_i$, represents the kinetic energy of the $i$-th particle, and $W\colon\R^d\to[0,\infty]$ is the trap potential. The trace of the operator ${\rm e}^{-\beta \Hcal_N}$ is the canonical partition sum of the system at temperature $1/\beta$. However, the characteristic property of bosons is expressed by the symmetry of any $N$-particle wave function under permutation of the coordinates. This in turn means that the partition sum of a system of $N$ bosons is given by the trace of the restriction of ${\rm e}^{-\beta \Hcal_N}$ to the subspace of symmetric wave functions, denoted by ${\rm Tr}_+({\rm e}^{-\beta \Hcal_N})$. Via the Feynman-Kac formula, this trace is given as
\begin{equation}\label{trace}
\begin{aligned}
{\rm Tr}_+\big({\rm e}^{-\beta \Hcal_N}\big)&=\frac 1{N!}\sum_{\sigma\in\Sym_N}\int_{\R^{dN}}\d x_1\cdots \d x_N\int_{\Ccal^N}\Big(\bigotimes_{i=1}^N \mu^{\beta}_{x_i,x_{\sigma(i)}}\Big)(\d \omega)\,\exp\Big\{-\sum_{i=1}^N\int_0^\beta W(\omega_s^{\ssup i})\,\d s\Big\},
\end{aligned}
\end{equation}
where $\mu^{\beta}_{x,y}$ is the canonical Brownian bridge measure with initial site $x\in \R ^d$ and terminal site $y\in \R ^d$ (see \eqref{nnBBM} for a detailed definition of the Brownian bridge measure), and we wrote $\omega=(\omega^{\ssup 1},\dots,\omega^{\ssup N})$.

\noindent We therefore consider a system of $N$ Brownian motions under symmetrized conditions, and we are interested in rigorously describing the large $N$ behavior of the system in the presence of a trap potential. In this paper, we show that the averaged path behavior (in the sense of empirical measures) of the Brownian motions can be described by a simple stochastic diffusion. 

\noindent We also show that it solves a large deviation problem introduced by Schr\"odinger, which involves minimization of relative entropy with given marginals. This problem is known in the literature as the Schr\"odinger bridge problem. Let us give an overview of this concept: Fix a time horizon $\beta>0$, $\nu_1,\nu_2$ probability measures in $\R ^d$ and a Markov process $(X_t)_{t\ge 0}$ with underlying probability measure $\P$. The Schr\"odinger bridge problem consists of finding a measure $Q$ on the space of all paths $\omega:[0,\beta]\to \R ^d$ such that
\begin{equation}\label{SchrProblem}
H(Q|\P)=\min_{Q_0\sim \nu_1,Q_\beta\sim \nu_2}H(Q|\P),    
\end{equation}
where $Q_0$ and $Q_\beta$ are the distributions of the initial and final time of Q respectively and $H(Q|\P)$ is the relative entropy of $Q$ with respect to $\P$ (see \eqref{Entropydef}.) The authors of [FG97] have expressed the solutions of the Schr\"odinger bridge problem as an integral of bridge measures: Let $\P_{x,y}$ be the conditional probability measure $\P$, given that $X_0=x$ and $X_\beta=0$. Then we can write $\P$ as 
$$\P= \int_{\R ^d}\int_{\R ^d}\P_{x,y}\mu(\d x,\d y),$$
where $\mu$ is the joint law of $X_1$ and $X_2$ under $\P$. Let now $\nu $ be a probability measure on $\R ^d \times \R ^d$ with marginal distributions $\nu(\d x, \R ^d)=\nu_1(\d x)$ and $\nu(\R ^d,\d y)=\nu_2(\d y)$. If $(X_t)_{t\ge 0}$ under the measure 
$$Q:=\int_{\R ^d}\int_{\R ^d}\P_{x,y}\nu(\d x,\d y)$$
remains Markovian, then $Q$ is called a Schr\"odinger process. According to \cite{FG97} a Schr\"odinger process is a solution to the problem \eqref{SchrProblem}. In the Appendix~$\ref{sec-sch}$ we provide a detailed definition of a Schr\"odinger process.

\noindent Finally, we extend the Large deviations results of \cite{AK08} and we provide an explicit formula for the rate function that controls the mean of occupation measures defined in \eqref{YNdef}.

\noindent The paper is organized as follows: In section 2 we set up the model and state some auxiliary results, in section 3 we state our main results, in section 4 we discuss the connection of our problem with the existing literature and in section 5 we provide proofs of our main results. The appendix consists of proofs of the auxiliary results and detailed definitions of Large Deviations and Schr\"odinger processes. 

\section{Set up and auxiliary results}\label{sec- set up and results}\smallskip

\noindent We study the large-$N$ behavior of a system of $N$ symmetrized trapped Brownian motions $B_1,\dots, B_N$ in $\R^d$ on a fixed time interval $[0,\beta]$, i.e., the behavior of the system under the measure defined as follows: For a tuple of $ N $ paths $ \omega=(\omega^{\ssup{1}},\ldots,\omega^{\ssup{N}}) $ with $ \omega^{\ssup{i}} \in\Ccal_\beta:= C([0,\beta];\R^d) $ where $C([0,\beta];\R^d)$ is the space of all continuous paths $\omega:[0,\beta]\to \R ^d$, we set
\begin{equation}\label{transsym}
     \P_{\m,N,W}^{\ssup{{\rm sym}}}(\d \omega):=\frac{1}{Z_{\m,N}^{\ssup{{\rm sym}}}}\exp\left(-H_N(\omega)\right)  \P_{\m,N}^{\ssup{{\rm sym}}}(\d\omega), 
\end{equation}
where 
\begin{equation}\label{Psymdef}
\P_{\m,N}^{\ssup{{\rm sym}}}=\frac{1}{N!}\sum_{\s\in\Sym_N}\;\int_{\R^d}\cdots \int_{\R^d} \m(\d x_1)\cdots\m(\d x_N)\bigotimes_{i=1}^N\P^{\beta}_{x_i,x_{\s(i)}}, 
\end{equation}
\begin{equation}
    H_N(\omega):=\sum_{i=1}^N\int_0^\beta W(\omega^{(i)}(s))\d s
\end{equation}
and $Z_{\m,N}^{\ssup{{\rm sym}}}$ is the partition function defined as 
\[\Z_{\m,N}^{\ssup{{\rm sym}}}:=\E_{\P_{\m,N}^{\ssup{{\rm sym}}}}[e^{-H_N}].\]
Here $\Sym_N$ is the set of all permutations of $1,\dots,N$, $\P_{x,y}^{\beta}$ the normalised Brownian bridge measure on the time interval $[0,\beta]$ with initial site $ x \in\R^d$ and terminal site $ y \in\R^d$. Recall that the canonical (non-normalised) Brownian bridge measure on the time interval $[0,\beta]$ with initial site $x\in\R^d$ and terminal site $y\in \R^d$ is defined as
\begin{equation}\label{nnBBM}
\mu_{x,y}^\beta(A)=\frac {\P_x(B\in A;B_\beta\in\d y)}{\d y},\qquad A\subseteq \Ccal_\beta\mbox{ measurable,}
\end{equation}
where $\P_x$ is the probability measure under which $B=B^{\ssup{1}}$ starts from $x\in\R^d$. Then $\mu_{x,y}^\beta$ is a regular Borel measure on $\Ccal_\beta$ with total mass equal to the Gaussian density,
\begin{equation}\label{Gaussian}
\mu_{x,y}^\beta(\Ccal_\beta)=p_\beta(x,y)=\frac {\P_x(B_\beta\in\d y)}{\d y}=(4\pi\beta)^{-d/2}{\rm e}^{-\frac 1{4\beta}|x-y|^2}.
\end{equation}
The normalized Brownian bridge measure is defined as $\P_{x,y}^\beta=\mu_{x,y}^\beta/p_\beta(x,y)$, which is a probability measure on $\Ccal_\beta$. Moreover, $\m$ is the initial probability distribution on $\R^d$. 

\noindent Hence, the terminal location of the $i$-th motion is affixed to the initial location of the $\sigma(i)$-th motion, where $\sigma$ is a uniformly distributed random permutation. Any of the $N$ paths is a Brownian motion with initial distribution $\m$, but with a peculiar terminal distribution at time $\beta$. We can conceive $\P_{\m,N}^{\ssup{{\rm sym}}}$ as a two-step random mechanism: First we pick a uniform random permutation $\sigma$, then we pick $N$ Brownian motions with initial distribution $\m$, and the $i$-th motion is conditioned to terminate at the initial point of the $\sigma(i)$-th motion, for any $i$.

\noindent The function $W:\R ^ d \to [0,\infty]$ is the so-called "trap potential" and we assume that 
\begin{equation}\label{trapcond}
\lim_{R\rightarrow\infty}\inf_{|x|>R}W(x)=\infty,    
\end{equation}
it is bounded from below, and it is continuous on the set $\{x\in \R ^ d: W(x)<\infty\}$ which is either $\R ^ d$, or a connected compact subset of $\R ^ d $.

\noindent We are interested in the large-$N$ behavior of the empirical path measure defined as follows: For a tuple of $ N $ paths $ \omega=(\omega^{\ssup{1}},\ldots,\omega^{\ssup{N}}) $, we define $L_N:\Ccal_\beta^N \to \Mcal_1(\Ccal_\beta)$ as
\begin{equation}\label{LNdef}
L_N(\omega)=\frac{1}{N}\sum_{i=1}^N\delta_{\omega^{\ssup{i}}}\in\Mcal_1(\Ccal_\beta),
\end{equation}
where for any topological space $X$ we denote by $\Mcal_1(X)$ the set of all probability measures defined on $X$ and for $x \in X$ we denote by $\delta_x \in \Mcal_1(X)$ the point-mass probability measure.

\noindent More precisely, we show that the distributions of $L_N$ under $\P_{W,\m,N}^{\ssup{{\rm sym}}}$ as $N\to\infty$, converge in probability to a unique measure $\mu^*_{W,\m}$, and under further assumptions on the initial distribution $\m$, we show that $\mu^*_{W,\m}$ belongs to a class of measures that solve the so-called Schr\"odinger bridge optimal transport problem: an entropy minimization problem on path measures with fixed initial and terminal distributions (see Section~\ref{sec-scr} for a detailed introduction to Schr\"odinger processes and  the Schr\"odinger bridge problem). To reach this goal we extend the Large Deviations results in \cite{AK08}. In Section~\ref{sec-LDP}, we recall the notion of a Large Deviations Principle. In particular, we obtain a Large Deviations Principle for $L_N$ under the measure $\P_{W,\m,N}^{\ssup{{\rm sym}}}$ and for the means of the normalized occupation measures, defined as follows: For a tuple of $ N $ paths $ \omega=(\omega^{\ssup{1}},\ldots,\omega^{\ssup{N}}) $ we define $Y_N:\Ccal_\beta^N \to \Mcal_1(\R ^d )$ as
\begin{equation}\label{YNdef}
Y_N(\omega)= \frac{1}{N}\sum_{i=1}^N\frac{1}{\beta}\int_0^{\beta}\d s\,\delta_{\omega^{\ssup{i}}(s)}\in\Mcal_1(\R^d),
\end{equation}

\noindent We also study an important special case: We set $\m$ to be the Lebesgue measure and we replace $\P_{x,y}^\beta$ by the canonical, non-normalized, Brownian bridge measure, $\mu_{x,y}^\beta$ (see \eqref{nnBBM}). In this case, the rate function for the means of occupation measures $Y_N$ turns out to be equal to $ \beta $ times the well-known Donsker-Varadhan rate function, which we define in Section~\ref{sec-LDP}. 
\smallskip\\

The first step toward our results is to study Large Deviations Principles for the distributions of $L_N$ and $ Y_N $ under the measure defined in \eqref{transsym}. Throughout the paper, we fix $\beta >0$; we endow $\Ccal_\beta$ with the topology of uniform convergence and the corresponding Borel $\sigma$-algebra. We consider $N$ random variables, $ B^{\ssup{1}},\ldots, B^{\ssup{N}}$, taking values in $\Ccal_\beta$. For completeness, since we extend the results in \cite{AK08}, we work with Brownian motions having generator $\Delta$ instead of $\frac 12\Delta$. Moreover, for a function $f:\R ^ d\to \R$, we define 
\[\Tilde{f}(\omega):=\int_0^\beta f(\omega(s)) \d s,\qquad \omega \in \Ccal_\beta.\]
Additionally, we endow $\Mcal_1(\Ccal_\beta)$ with the L\'evy metric defined as 
\begin{equation}\label{levym}
    \d(\mu,\nu):=\inf\{\delta >0\colon\mu(\Gamma)\le \nu(\Gamma^{\delta})+\delta\,\mbox{ and }\, \nu(\Gamma)\le\mu(\Gamma^{\delta})+\delta\,\mbox{ for all }\, \Gamma=\overline{\Gamma}\subset\Ccal_\beta\}, 
\end{equation} 
where $F^{\delta}=\{\omega \in \Ccal_\beta \colon\dist(\omega,F)\leq\delta\} $ (here, $\dist(\cdot,\cdot)$ denotes the distance function that induces uniform convergence) is the closed $\delta$-neighborhood of $ F $. According to \cite{B68}, a sequence $(\mu_N)_{N\in \N}$ of elements in $\Mcal_1(\Ccal_\beta)$ converges in the d-metric to a measure $\mu$ if and only if it converges in the weak sense, i.e., for every $F \in C_b(\Ccal_\beta)$ (the space of all continuous and bounded functions on $\Ccal_\beta$) we have that
\[\int_{\Ccal_\beta}F\d \mu_N\to\int_{\Ccal_\beta}F\d \mu \]
as $N\to \infty$. Hence, the L\'evy metric generates the weak topology of probability measures. The latter is the topology under which we develop Large Deviations Principles. We also use the notation "$\mu_N\Rightarrow \mu$" to denote the weak convergence of probability measures.

\noindent Finally, let $X$ be a topological space. With 
\begin{equation}\label{Entropydef}
H(q | \widetilde{q})=\int_{X}q(\d x)\log\frac{q(\d x)}{\widetilde{q}(\d x)}q(\d x)
\end{equation} 
we denote the relative entropy of $q\in \Mcal_1(X)$ with respect to $ \widetilde q\in\Mcal_1(X)$. We often use this notation for $X=\R^d\times\R^d$ in the sequel, but also for other spaces $X$

\noindent We now recall the rate functions introduced in \cite{AK08}. Let  $ \Mcal_1^{\ssup{\rm s}}(\R^d\times\R^d)$ be the set of shift-invariant probability measures $q$ on $ \R^d\times\R^d $, i.e., measures whose first and second marginals coincide and are both denoted by $\overline q$. Note that $q\mapsto  H(q|\overline{q}\otimes\m)$ is strictly convex. We write $\langle F ,\mu\rangle$ for integrals $\int_{\Ccal_\beta} F(\omega)\,\mu(\d \omega)$ for suitable functions $F$ on $\Ccal_\beta$. 

\noindent Define the functional $I^{\ssup{{\rm sym}}}_{\m}$ on $ \Mcal_1(\Ccal_\beta) $ by
\begin{equation}\label{Isymdef}
I^{\ssup{{\rm sym}}}_{\m}(\mu)=\inf_{q\in\Mcal_1^{\ssup{\rm s}}(\R^d\times\R^d)}\Bigl\{H(q|\overline{q}\otimes\m)+ I^{\ssup{q}}(\mu)\Bigr\},\qquad \mu\in\Mcal_1(\Ccal_\beta),
\end{equation} 
where
\begin{equation}\label{Jqdef}
I^{\ssup{q}}(\mu)=\sup_{F \in\Ccal_{\rm b}(\Ccal)}\Bigl\{\langle F ,\mu\rangle -\int_{\R^d}\int_{\R^d} \log\E_{x,y}^{\beta}\big[{\rm e}^{F }\big]q(\d x, \d y)\Bigr\},\qquad\mu\in\Mcal_1(\Ccal_\beta).
\end{equation}
\begin{prop}[Large deviations for $L_N$]\label{main}
Fix $ \beta \in (0,\infty) $. Then, under the assumption \eqref{trapcond} (and the subsequent assumptions)  as $ N\to\infty$, under the measure $ \P_{\m,N,W}^{\ssup{{\rm sym}}} $ defined in \eqref{transsym}, the empirical path measures $ L_N $ satisfy a Large Deviations Principle on $ \Mcal_1(\Ccal_\beta) $ with speed $ N $ and rate function $I_{W,\m}:\Mcal_1(\Ccal_\beta)\to \R \cup \{\infty\}$ given by
\begin{equation}\label{Trapratefuncemp}
I^{(sym)}_{\m,W}(\mu):=I^{(sym)}_{\m}(\mu)+\langle \Tilde{W},\mu\rangle-\inf_{\mu \in\Mcal_1(\Ccal_\beta)}\Bigl\{I^{(sym)}_{\m}(\mu)+\langle \Tilde{W},\mu\rangle\Bigr\}.
\end{equation}
\end{prop}
\noindent To be more explicit, the stated Large Deviations Principle says that
$$
\lim_{N\to\infty}\frac 1N\log \P_{\m,N,W}^{\ssup{{\rm sym}}}\big(L_N\in\,\cdot\big)=-\inf_{\mu\in\,\cdot}I^{\ssup{{\rm sym}}}_{W,{\m}}(\mu),
$$
in the weak sense, i.e., there is a lower bound for open subsets of $\Mcal_1(\Ccal_\beta)$ and an upper bound for closed ones. The proof of Proposition~\ref{main} is in the appendix~\ref{sec-auxresults}. This extends Theorem 1.1 in \cite{AK08}. Note that in the latter, it is assumed that the initial distribution $ \m\in\Mcal_1(\R^d) $ has compact support. We show that this assumption can be relaxed and $\m$ can also be supported in $\R ^d$. 

\noindent We also have an analogous argument for the mean of the occupation measures, $Y_N$, defined in \eqref{YNdef}.
To formulate this, we define the functional $J^{\ssup{\rm sym}}_{\m}$ on $ \Mcal_1(\R^d) $ by
\begin{equation}\label{Jsymdef}
J^{\ssup{\rm sym}}_{\m}(p):=\inf_{q\in\Mcal_1^{\ssup{\rm s}}(\R^d\times\R^d)}\Bigl\{H(q|\overline{q}\otimes\m)+ J^{\ssup{q}}(p)\Bigr\},\qquad p\in\Mcal_1(\R^d),
\end{equation}
where
\begin{equation}\label{Jsymqdef}
J^{\ssup{q}}(p):=\sup_{f\in\Ccal_{\rm b}(\R^d)}\Bigl\{\beta\langle f,p\rangle -\int_{\R^d}\int_{\R^d} q(\d x,\d y)\,\log\E_{x,y}^{\beta}\Bigl[{\rm e}^{\int_0^{\beta}f(B_s)\,\d s}\Bigr]\Bigr\}.
\end{equation}

\begin{cor}[Large deviations for $Y_N$]\label{main2}
Fix $ \beta \in (0,\infty) $. Then, as $ N\to\infty $, under the measure $ \P_{\m,N,W}^{\ssup{{\rm sym}}} $ defined in \eqref{transsym}, the mean of occupation measures, $ Y_N $, satisfy a Large Deviations Principle on $ \Mcal_1(\R^d) $ with speed $ N $  and rate function $ J^{\ssup{\rm sym}}_{W,\m} :\Mcal_1(\R^d)\to \R \cup \{\infty\}$ given by:
\begin{equation}\label{Trapratefuncocp}
J^{(sym)}_{\m,W}(p):=J^{(sym)}_{\m}(p)+\langle W,p\rangle-\inf_{p\in\Mcal_1(\R^d)}\Bigl\{J^{(sym)}_{\m}(p)+\langle W,p\rangle\Bigr\}.
\end{equation}

\end{cor}

\noindent The proof of Corollary~\ref{main2} can be found in the appendix~\ref{sec-auxresults}. In \cite{AK08} the initial probability measure $\m \in \Mcal_1(\R ^ d)$ is extended to any finite measure in $\R ^d$. Moreover, the symmetrized measure $\P^{\ssup{\rm sym}}_{\m,N}$ is also extended to a larger class of measures. We can use similar arguments to extend the results of Proposition~\ref{main}, and Corollary~\ref{main2}. Note that here, $\m$ is extended to a larger class of measures.   
\begin{prop}\label{cormain}
Fix $ \beta \in (0,\infty) $ and assume that $ \m $ is either a positive finite measure on $\R^d$, or a $\sigma-$ finite measure which is absolutely continuous with respect to the Lebesgue measure, with continuous and bounded density. Fix some continuous and bounded function $g\colon\R^d\times \R^d\to(0,\infty)$ and replace $\P_{x,y}^\beta$ by $g(x,y)\P_{x,y}^\beta$ in the definition \eqref{Psymdef} of $\P_{\m,N}^{\ssup{{\rm sym}}}$. Then 
\begin{enumerate}
\item[(i)] Proposition~\ref{main} remains true. The corresponding rate function is
\begin{equation}\label{Isymgdef}
    \mu\mapsto I^{\ssup{\rm sym}}_{W,\m,g}(\mu)-\langle \mu_{0,\beta},\log g\rangle-\inf_{\mu \in \Mcal_1(\Ccal_\beta)}\Bigl\{I^{\ssup{\rm sym}}_{W,\m}(\mu)-\langle \mu_{0,\beta},\log g\rangle\Bigr\}.
\end{equation}

\item[(ii)] Corollary~\ref{main2} remains true. The corresponding rate function is
\begin{equation}\label{Jsymgdef}
J^{\ssup{\rm sym}}_{\m,W,g}(p)=\widetilde{J}_{g,W,N}(p)-\inf_{p\in \Mcal_1(\R ^ d)}\{\widetilde{J}_{g,W,N}(p)\},\qquad p\in\Mcal_1(\R^d),
\end{equation}
where 
\[\widetilde{J}_{g,W,N}(p)=\inf\limits_{q\in\Mcal_1^{\ssup{{\rm s}}}(\R^d\times\R^d)}\Big\{H(q|\overline{q}\otimes\m)+J^{\ssup{q}}(p)+\langle W, p\rangle-\langle q,\log g\rangle\Big\}.\]
\end{enumerate}
\end{prop}

\section{Results}\label{sec-results}\bigskip

\noindent We now  proceed to our main results. We first focus on the rate function introduced in part (i) of Proposition \ref{cormain}. In particular, we show that under the assumptions of Proposition \ref{cormain} and for any choice of the function $g\in C_b(\R^d\times\R^d)$ and the initial measure $\m\in \Mcal(\R^d) $ introduced in Proposition \ref{cormain}, the mapping $\Mcal_1(\Ccal_\beta)\ni\mu \mapsto I^{\ssup{\rm sym}}_{\m,W,g}(\mu)$ has a unique minimizer. 
\noindent We state our results. Here the trap $W$ is any function that satisfies the assumptions in Theorem \ref{main}. Here we only consider the non-normalized rate function i.e., $I^{(sym)}_{W,\m}(\mu):=I^{(sym)}_{\m}(\mu)+\langle \Tilde{W},\mu\rangle$. A key part of our results is the following:
\begin{prop}\label{unmin}
Under the assumptions of Proposition \ref{cormain}, the following hold
\begin{enumerate}

\item[(i)] We have that
\begin{equation}\label{infIg}
\begin{aligned}
&\inf_{\mu \in \Mcal_1(\Ccal_\beta)}\Bigl\{I^{\ssup{\rm sym}}_{W,\m}(\mu)-\langle \mu_{0,\beta},\log g\rangle\Bigr\}=\\
&\inf\limits_{q\in\Mcal_1^{\ssup{{\rm s}}}(\R^d\times\R^d)}\Big\{H(q|\overline{q}\otimes\m)-\int_{\R ^d}\int_{\R ^d}\log \E^\beta_{x,y}\Big[e^{-\Tilde{W} }\Big]q(\d x, \d y)-\langle g,\log g\rangle\Bigr\}.\\
\end{aligned}
\end{equation}

\item[(ii)] The unique minimizer of the (non-normalized) rate function defined in Proposition \ref{cormain} (i) is given by 
\begin{equation}\label{minimizergen}
    \mu^*_{g,\m,W}=\int_{\R^d}\int_{\R^d}Q^{\beta}_{x,y}\circ B^{-1}q^*_{g,\m,W}(\d x, \d y),
\end{equation}
where for $x,y \in \R ^d $, we define
\begin{equation}\label{Qmeasure}
    Q^\beta_{x,y}(\d \omega):=\frac{1}{\E^\beta_{x,y}\big[e^{-\Tilde{W} }\big]}e^{-\Tilde{W}(\omega)}\P^\beta_{x,y}(\d \omega),\qquad \omega \in \Ccal_\beta,
\end{equation}
and $q_{g,\m,W}^*\in\Mcal_1^{\ssup{{\rm s}}}(\R^d\times\R^d)$ is the unique minimizer of the formula on the right hand side of \eqref{infIg}.
\end{enumerate}
\end{prop}
\noindent The proof of Proposition~\ref{unmin} can be found in Section~\ref{sec-mintheorem}. Additionally, if we impose further conditions on $g,\mu$ and $W$, we can simplify the results of Theorem~\ref{unmin}:
\begin{cor}\label{eigencor}
    Assume in addition that $g,\mu$ and $W$ are such that the operator 
    \begin{equation}\label{Tdef}
        T(f)(x):=\int_{\R ^d} f(y) \E_{x,y}^{\beta}\left [e^{-\Tilde{W} }\right]g(x,y)\m(\d y)
    \end{equation}
    attains its minimum eigenvalue $\lambda_T>0$ with corresponding $\m$-normalized eigenfunction $\varphi^2_T \in H^1(\R ^d)$. Then \item[(i)] we have
\begin{equation}
\begin{aligned}\label{infIgT}
&\inf_{\mu \in \Mcal_1(\Ccal_\beta)}\Bigl\{I^{\ssup{\rm sym}}_{W,\m}(\mu)-\langle \mu_{0,\beta},\log g\rangle\Bigr\}=-\log \lambda_T. 
\end{aligned}
\end{equation}
\item[(ii)] The minimizer of (i) is given by 
\begin{equation}\label{minT}
    \mu^*_{T}=\frac{1}{\lambda_T}\int_{\R^d}\int_{\R^d} \varphi_T(x)\varphi_T(y)e^{-\Tilde{W} }\P_{x,y}^{\beta}g(x,y)\m(\d y)\m(\d x).
\end{equation}
\end{cor}
\smallskip
\noindent For some special cases of the function $g$, the initial measure $\m$ and the trap potential $W$, the minimizer $\mu^*_{g,\m,W}$ is identified with the solution of the Schr\"odinger Bridge problem as formulated in [FG97]. In particular, if we set $\m=\Leb$, the Lebesgue measure, $g\equiv p_\beta$ the Gaussian density, and $W$ to be either a Hard-wall or Soft-wall potential then $Q_{x,y}$ in \eqref{Qdef} depends only on the initial and terminal projections of the path. \textbf{Hard-wall potentials} are potentials of the form $W=\1_{\L^c}\infty$, where $\L$ is a compact box in $\R^d$, and \textbf{Soft-wall potentials} are functions such that $W<\infty$ everywhere, and if $\varphi$ is the eigenfunction of $-\Delta+W$ associated with its minimum eigenvalue, then we have that $W\varphi \in L^2(\R^d)$. According to [RS78IV], Theorem XIII.72, if $W>0$, and smooth, then the latter assumption is satisfied. This property ensures that $\varphi \in H^2(\R ^d)$. Our first result is a simplified expression of $\mu^*_{g,\m,W}$

\begin{theorem}\label{propmininmizer}
    Let $\L$ be a bounded closed box, and let $W$ be a Hard-wall potential. Then the measure $\mu^*_{g,\m,W}$ defined in \eqref{minimizergen} is an ergodic Markov process with invariant measure being the uniform distribution on $\L$. Moreover, it is the unique minimizer of the rate function $I_{\m,W,g}^{(sym)}$, when $\m=\Leb_\L$ and $g\equiv p_\beta$. Finally, $\mu^*_{g,\m,W}$ solves the optimal transport problem $$\inf\{H(Q|\P_\m): Q\circ\pi_0^{-1}=U_{\L},Q\circ\pi_\beta^{-1}=U_{\L}\},$$
Where $U_{\L}$ denotes the Uniform distribution on $\L$
\end{theorem}

\noindent The proof of Theorem~\ref{propmininmizer} can be found in section~\ref{sec-proof}. This result aligns with the results in \cite{AK08}. The authors showed that, in the case of a hard-wall potential, the rate function \eqref{Jsymgdef} is given by $J(p)=\beta\|\nabla \phi\|_2^2$, when $\frac{\d p}{\d x}=\phi^2\in H^1_0(\L^o)$ and $J(p)=\infty$ otherwise. The unique minimizer of this rate function is $\phi^2(x)=\1_{\L^o}(x)\frac{1}{|\L|}$ which is the density of the uniform distribution on $\L$. 

\noindent We also have the analogous result for soft-wall potentials. Since in this case, the trap potentials belong to a broader class of functions, the following result becomes less explicit than the one in Theorem~\ref{propmininmizer}.
\begin{theorem}\label{mintheorem}
    Fix $\beta \in (0,\infty)$, and let $\tilde W:\Ccal_\beta\to \R$ defined as $\tilde W(\omega)=\int_0^\beta W(\omega(s)) \d s$ for all $\omega \in \Ccal_\beta$, where $W$ is a Soft-wall potential. Denote by  $\varphi:\R ^ d \to \R$ the $L^2(\R^d)$-normalized eigenfunction of the operator $-\Delta +W$ which corresponds to the principle eigenvalue 
        \begin{equation}\label{eigenvdefi}
        \l(W)=\sup_{\varphi\in H^1(\R^d)\colon , \|\varphi\|_2=1}\Big(\langle W,\varphi^2\rangle-\|\nabla\varphi\|_2^2\Big).
    \end{equation}
    \noindent Then, setting
    \begin{equation}\label{pairscr}
     q^*(\d x,\d y):=\varphi(x)\varphi(y)\frac{\E_x\left[e^{\tilde W -\beta\l (W)}; B_\beta \in \d y\right]}{\d y}\,\d x \d y,  
    \end{equation}
 we have that $\mu^*_{g,\m,W}$ defined in \eqref{minimizergen} with $q^*_{g,\m,W}$ replaced by $q^*$ in \eqref{pairscr} is an ergodic Markov process which corresponds to the solution to the problem
 \begin{equation}\label{sdesh}
     \d X_t=\frac{\nabla \varphi}{\varphi}\d t+\d B_t,\quad X_0 \text{ has distibution with density } \varphi^2 \text{ with respect to the Lebesgue measure}
 \end{equation}
 with invariant measure with density $\phi^2$, and it is the unique minimizer of $I_{\m,W,g}^{(sym)}$, when $\m=\Leb$ and $g\equiv p_\beta$. Moreover, $\mu^*_{g,\m,W}$ is a solution to the optimal transport problem:
$$\inf\{H(Q|\P_\m): Q\circ\pi_0^{-1}=\phi^2(x)\d x,Q\circ\pi_\beta^{-1}=\phi^2(x)\d x\}.$$
\end{theorem}
\noindent The proof of Theorem~\ref{mintheorem} can be found in section \ref{sec-sch}.

\noindent Note that the principal eigenvalue \eqref{eigenvdefi} is the minimum value of the rate function that controls the mean of empirical measures $L_N$, and therefore it is the minimum value of the rate function that controls the mean of occupation measures $Y_N$. The expression in the supremum in \eqref{eigenvdefi} when multiplied by $\beta$ turns out to be the rate function the controls $Y_N$. This form of rate function is called the Donsker-Varadhan rate function.  
\subsection{An important special case}\bigskip

\noindent Consider the non-normalized rate function 
\begin{equation}\label{Jnonnormal}
   \widetilde{J}_{W,N}(p)=\inf\limits_{q\in\Mcal_1^{\ssup{{\rm s}}}(\R^d\times\R^d)}\Big\{H(q|\overline{q}\otimes\m)+J^{\ssup{q}}(p)+\langle W, p\rangle-\langle q,\log g\rangle\Big\}.
\end{equation}
Here the trap potential $W$ is a Soft-wall potential. We choose the initial measure to be $\m=\Leb$ the Lebesgue measure on $\R^d$ and $g\equiv p_\beta$ the Gaussian density defined in \eqref{Gaussian}. Denote by $\overline{J}$ the rate function \eqref{Jnonnormal} when we choose $\m=\Leb$ and $g\equiv p_\beta$. Our result is an extension of Theorem 1.5 in \cite{AK08}.
\begin{theorem}\label{thDVtrap}
    Consider the mapping $I_W:\Mcal_1(\R^d)\to [0,\infty]$,
    \begin{equation}\label{DVtrap}
    I_W(p):=
        \begin{cases}
            \|\nabla \phi\|_2^2+\langle W,\phi^2\rangle, \text{ if } \frac{\d p}{\d x}=\phi^2 \in H^1(\R^d)\\
            \infty, \text{ otherwise.}
        \end{cases}
    \end{equation}
    Then, $\overline{J}(p)=\beta I_W(p)$
\end{theorem}
\noindent The proof of Theorem~\ref{DVtrap} can be found in Section \ref{sec-DV}

\section{Relation to quantum statistical mechanics}\label{sec-Quantum}\bigskip

\noindent In the Introduction of this paper, we described the relation of our work with the canonical ensemble of large systems of bosons at positive temperature. Ultimately, this project is a stepping stone toward a complete probabilistic interpretation of Bose-Einstein condensation. Let us recall some of the existing results in the literature and some of our future aims. In \cite{AK08}, the authors proved an explicit formula for the logarithmic large-$N$ asymptotic of the trace formula \eqref{trace}, for a certain class of hard-wall traps $W$, in particular, when $W=\1_{\L^c}\infty$, where $\L$ is a finite box of $\R ^ d$: The large-$N$ behavior of $L_N$ defined in \eqref{LNdef} and $Y_N$ defined in \eqref{YNdef} has been studied, under the measure defined in \eqref{Psymdef}. In particular, if the initial measure $\m\in \Mcal_1(\R ^ d)$ has compact support, then $L_N$ satisfies a Large Deviations Principle with speed $N$ and rate function $I^{\ssup{{\rm sym}}}_{\m}$ defined in \eqref{Isymdef}, and $Y_N$ satisfies a Large Deviations Principle with speed $N$ and rate function $J^{\ssup{{\rm sym}}}_{\m}$ defined in \eqref{Jsymdef}. Moreover, under the assumptions of Proposition \ref{cormain}, and if we assume that the initial measure $\m\in \Mcal(\R ^ d)$ has compact support, then $L_N$ satisfies a Large Deviations Principle with speed $N$ and rate function given by the mapping $ \mu\mapsto I^{\ssup{\rm sym}}_{\m}(\mu)-\langle \mu_{0,\beta},\log g\rangle$ $\mu \in \Mcal_1(\Ccal_\beta)$, and $Y_N$ satisfies a Large Deviations Principle with speed $N$ and rate function given by the mapping $p \mapsto \inf\limits_{q\in\Mcal_1^{\ssup{{\rm s}}}(\R^d\times\R^d)}\Big\{H(q|\overline{q}\otimes\m)+J^{\ssup{q}}(p)-\langle q,\log g\rangle\Big\},\quad p\in\Mcal_1(\R^d)$. A special case related to the partition function \eqref{trace} is when we set $\m=\Leb_{\L}$, where $\L \subset \R ^d$ is a connected compact subset of $\R ^d$, and $g\equiv p_\beta$, where $p_\beta$ is the Gaussian density function (see \eqref{Gaussian}). Then the rate function corresponding to the Large Deviations Principle satisfied by $Y_N$ under the measure defined in proposition \ref{cormain} is $\beta$ times a functional that belongs to well known class of functionals called Donsker-Varadhan rate functions (see Section \ref{sec-LDP}), more precisely if we define the functional $I_\L\colon \Mcal_1(\R^d)\to[0,\infty]$ by
\begin{equation}\label{Idef}
I_\L(p)=\begin{cases}\big\|\nabla \sqrt{\frac{\d p}{\d x}}\big\|_2^2,&\mbox{if }p\mbox{ has a density with square root in }H^1_0(\L^{\circ}),\\
\infty&\mbox{otherwise,}
\end{cases}
\end{equation}
then, the corresponding rate function is $p\mapsto \beta I_\L(p)$. Finally, since $W=\1_{\L ^c}$, then $W$ is continuous and bounded in $\L$, thus $\Mcal_1(\L)\ni Y\mapsto \langle W,Y\rangle$ is bounded and continuous, we may apply Varadhan's lemma to deduce that
$$
\lim_{N\to\infty}\frac 1N \log\Big({\rm Tr}_+\big({\rm e}^{-\beta \Hcal_N}\big)\Big)=-\inf_{p\in \Mcal_1(\L)}\Big(J^{\ssup{\rm sym}}_{\L}(p)+\beta\langle W,p\rangle\Big)=\beta\lambda_\L(W),
$$
where $\lambda_\L(W)$ is defined in \eqref{eigenvdefi}. This result is extended in the case when $W$ is a soft wall trap potential. 

\noindent We also aim to extend our results by adding to the Hamiltonian operator in \eqref{Hdef} an interaction function $ v\colon[0,\infty)\to\R $ that satisfies some suitable regularity conditions. The Hamilton operator now reads
\begin{equation}\label{NewHamilt}
    H_N=\sum_{i=1}^N-\Delta_i+W(x_i)+\sum_{1\leq i<j\leq N}u(x_i-x_j).
\end{equation}
Moreover, we are interested in comparing our results with the ones in \cite{MU11} when $d=3$. In their work, the authors proved that we can interpret the stochastic mechanics of an $N$-body system of bosons as a single particle process in a Bose-Einstein condensate. Let us briefly recall their results. Consider the N-body Hamiltonian defined in \eqref{NewHamilt}, and denote by $\Psi_N$ its ground state which we assume that it is strictly positive and of class $C^1$. The corresponding $3N$-dimensional Nelson's diffusion is given by:
\[\d \hat{X}_t=\frac{\nabla^{(N)}\Psi_N}{\Psi_N}(\hat{X}_t)\d t+\d B_t,\]
where $\hat{X}=(X_1,\dots ,X_N)$, with $X_i$ being a single particle process and $\nabla^{(N)}$ denotes the $3N$-dimensional gradient. If Bose-Einstein condensation occurs, then then condensate is described by the order parameter $\phi_{GP}\in L^2(\R^3)$ which is the minimizer of
\begin{equation}\label{phigp}
    \inf_{\{\varphi\in L^2(\R^3)\colon \|\varphi\|_2^2=1\}}\{\|\nabla \varphi\|_2^2+\langle W,\varphi\rangle +\frac \alpha2 \|\varphi\|_4^4\},
\end{equation}
where $\alpha$ is a parameter that depends on $v$, and the functions $W$ and $v$ satisfy some suitable regularity conditions. Define now the following diffusion
\begin{equation}\label{BEcon}
\d (X_{GP})_t=u_{GP}((X_{GP})_t)\d t +dB_t,    
\end{equation}
where 
\[u_{GP}(x)=\frac{\nabla \phi_{GP}}{\phi_{GP}}(x).\]
The authors showed that the one particle conditional distribution $X_1$ converges to $X_{GP}$ as $N\to \infty$ in a relative entropy sense. We conjecture that that the empirical path measure $L_N$ under the transformed symmetrized measure with density $e^{-K_N-H_N}$, where $K_N(\omega):=\sum_{1\le i<j\le N} \frac{1}{\beta}   \int_0^\beta\int_0^\beta \, v\big(|\omega^{\ssup{i}}_s-\omega^{\ssup{j}}_t)|\big)\,\d s\d t$ and $H_N(\omega)= \sum_{i=1}^N \int_0^\beta\; W(\omega^{\ssup{i}}_s)\,\d s$ converges weakly to the solution of \eqref{BEcon}.\bigskip

\section{Proofs}\label{sec-proof}
\subsection{Proof of Proposition~\ref{mintheorem}.}\label{sec-mintheorem} We split the proof into two Lemmas and a corollary. In the first lemma we prove claim (i), in the second lemma we set up the proof of claim (ii) which we eventually prove in Corollary~\ref{corproofmin}.

\noindent Let $q\in \Mcal_1(\R ^d \times \R ^d)$. We use the notation
\begin{equation}\label{fltransform}
 I_W^\ssup{q}(\mu):=\sup_{F \in \Ccal_\beta}\left\{\langle F+\tilde W,\mu\rangle -\int_{\R ^d}\int_{\R ^d}\log\E_{x,y} \left[e^F\right]q(\d x , \d y)\right\}.    
\end{equation}

\begin{lemma}\label{minbd}
We have
    \begin{equation}\label{min}
        \inf_{\mu \in \Mcal_1(\Ccal_\beta)}\{I_{\m,W,g}^\ssup{sym}(\mu)\}=\inf_{q\in \Mcal^{(s)}_1(\R ^ d\times \R ^d)}\left\{H(q|\overline{q}\otimes \m) -\int_{\R^ d} \int_{\R ^d}\log \E_{x,y}\left[e^{-{\tilde W }}\right]q(\d x,\d y)-\langle \log g,q\rangle \right \}.
    \end{equation}
\end{lemma}
\begin{proofsect}{Proof} Following a similar reasoning as in ([AK08], corollary 1.4), for $\mu \in \Ccal_\beta$, and $q\in \Mcal^{(s)}_1(\R ^ d\times \R ^d)$, we set $q=\mu_{0,\beta}$ in \eqref{Isymgdef} (otherwise $I_{\m,W,g}^\ssup{sym}(\mu)=\infty$, see [AK08], eq. 1.9). Moreover, we have 
\[\inf_{\mu \in \Mcal_1(\Ccal_\beta)}\{I_{\m,W,g}^\ssup{sym}(\mu)\}=\]
\[\inf_{q\in \Mcal^{(s)}_1(\R ^ d\times \R ^d)}\left\{H(q|\overline{q}\otimes \m)+\inf_{\mu \in \Mcal_1(\Ccal_\beta),\mu_{0,\beta}=q}I_{W}^\ssup{q}(\mu)-\langle \log g,q\rangle\right\}.\]
To see the this, let $q\in \Mcal^{(s)}_1(\R ^ d\times \R ^d)$ and any $\mu \in \Mcal_1(\Ccal_\beta)$. Then, 
\[H(q|\overline{q}\otimes \m)+\inf_{\mu \in \Mcal_1(\Ccal_\beta),\mu_{0,\beta}=q}I_{W}^\ssup{q}(\mu)-\langle \log g,q\rangle\leq H(q|\overline{q}\otimes \m)+I_{W}^\ssup{q}(\mu)-\langle \log g,q\rangle.\]
The latter holds for all $q\in \Mcal^{(s)}_1(\R ^ d\times \R ^d)$ and any $\mu \in \Mcal_1(\Ccal_\beta)$, since if $\mu_{0,\beta }\neq q$ the R.H.S is infinity and the inequality holds trivially. Hence, we have that
\[\inf_{q\in \Mcal^{(s)}_1(\R ^ d\times \R ^d)}\big\{H(q|\overline{q}\otimes \m)+\inf_{\mu \in \Mcal_1(\Ccal_\beta),\mu_{0,\beta}=q}I_{W}^\ssup{q}(\mu)-\langle \log g,q\rangle\big\}\]
is a lower bound of the set $\{I_{\m,W,g}^\ssup{sym}(\mu):\mu \in \Mcal_1(\Ccal_\beta)\}$. Let now $\epsilon>0$. For $q\in \Mcal^{(s)}_1(\R ^ d\times \R ^d)$, pick $\mu_q$ such that 
\[I_{W}^\ssup{q}(\mu_q)< \inf_{\mu \in \Mcal_1(\Ccal_\beta),\mu_{0,\beta}=q}I_{W}^\ssup{q}(\mu)+\epsilon,\]
which implies that
\[\inf_{\mu \in \Mcal_1(\Ccal_\beta)}I_{\m,W}^\ssup{sym}(\mu)\leq H(q|\overline{q}\otimes \m)+\inf_{\mu \in \Mcal_1(\Ccal_\beta),\mu_{0,\beta}=q}I_{W}^\ssup{q}(\mu)-\langle \log g,q\rangle +\epsilon.\]
Since the latter holds for any $q\in \Mcal^{(s)}_1(\R ^ d\times \R ^d)$, the assertion follows. 

\noindent We now prove that 
\begin{equation}\label{minIq}
    \inf_{\mu \in \Mcal_1(\Ccal_\beta),\mu_{0,\beta}=q}I_W^\ssup{q}(\mu)=-\int_{\R^ d} \int_{\R ^d}\log \E_{x,y}\left[e^{-{\tilde W}}\right]q(\d x,\d y).
\end{equation}
According to Lemma~\ref{lemmaalternativetrap}, we can write $I_{\m,W}$ as $ \widetilde{I}^{(sym)}_{\m,W}$ as seen in \eqref{Itildesymdef}. We now see that the inequality ``$\geq$" in \ref{minIq} holds if we choose $F\equiv 0$. For the ``$\leq$" direction of the inequalty, we choose $\mu=\int_{\R^ d} \int_{\R ^d} \Q^\beta_{x,y}q(\d x , \d y)$,
 where for $x,y \in \R ^d$
 \begin{equation}\label{Qdef}
      Q^{\beta}_{x,y}(\d \omega)= e^{-\Tilde{W}(\omega)}\frac{1}{\E_{x,y}[e^{-\tilde W}]}  \P^\beta_{x,y}(\d \omega), \qquad \omega \in \Ccal_\beta,
 \end{equation}
and adjust \eqref{fltransform} as follows:
 \[\sup_{F \in \Ccal_\beta}\left\{\langle F+\tilde W ,\mu\rangle -\int_{\R ^d}\int_{\R ^d}\log\E_{x,y} [e^F]q(\d x , \d y)\right\}=\]
 \[=\sup_{F \in \Ccal_\beta}\left\{\langle F+\tilde W ,\mu\rangle -\int_{\R ^d}\int_{\R ^d}\log\frac{\E_{x,y}[e^{F+\tilde W }e^{-\tilde W }]}{\E_{x,y}[e^{-\tilde W}]} q(\d x , \d y)\right\}-\int_{\R^ d} \int_{\R ^d}\log \E_{x,y}[e^{-{\tilde W }}]q(\d x,\d y)=\]
 \[\sup_{F \in \Ccal_\beta}\left \{\int_{\R^ d} \int_{\R ^d}\E_{Q^{\beta}_{x,y}}[ F+\tilde W]q(\d x, \d y) -\int_{\R ^d}\int_{\R ^d}\log\E_{Q^{\beta}_{x,y}}[e^{F+\tilde W }] q(\d x , \d y)\right\}-\int_{\R^ d} \int_{\R ^d}\log \E_{x,y}[e^{-{\tilde W }}]q(\d x,\d y),\]
  where $\E_{Q^{\beta}_{x,y}}$ is the expectation with respect to $Q^{\beta}_{x,y}$. We now see that $"\leq"$ holds using Jensen's inequality.  
 \end{proofsect}
\qed
\begin{lemma}\label{propminbd}
    The unique minimizer of \eqref{minIq} is given by 
    \begin{equation}\label{minimizer}
        \mu_q^*=\int_{\R^ d} \int_{\R ^d} \Q^\beta_{x,y}q(\d x , \d y),\quad q\in \Mcal_1^{(s)}(\R ^ d\times \R ^d),
    \end{equation}
    where for $x,y \in \R ^d$, $\Q_{x,y}$ is defined in \eqref{Qdef}
\end{lemma}
\begin{proofsect}{Proof} We first note that $I_W^\ssup{q}$ is convex as a combination of a supremum of linear functions, thus it attains a minimizer. Let $\mu \in \Mcal_1(\Ccal_\beta) $ be such that 
\begin{equation}
    I_W^\ssup{q}(\mu)=-\int_{\R^ d} \int_{\R ^d}\log \E_{x,y}[e^{-{\tilde W}}]q(\d x,\d y).
\end{equation}
Then, we have that
\begin{equation}\label{variationaltrap}
   \sup_{F \in \Ccal_\beta}\left\{\langle F,\mu\rangle -\int_{\R ^d}\int_{\R ^d}\log\E_{x,y} [e^{F-\tilde W}]q(\d x , \d y)\right\}=-\int_{\R^ d} \int_{\R ^d}\log \E_{x,y}[e^{-{\tilde W}}]q(\d x,\d y).
\end{equation}
We now show that the function 
\[J_{q,\mu}(F):=\langle F,\mu\rangle -\int_{\R ^d}\int_{\R ^d}\log\E_{x,y} [e^{F-\tilde W}]q(\d x , \d y),\quad F \in \Ccal_\beta,\]
is strictly concave. To see this, fix $F,G \in \Ccal_\beta$ and $\lambda\in (0,1)$ and note that using H\"older's inequality, we have
\begin{equation}
  \begin{aligned}
  &J_{q,\mu}(\l F+(1-\l)G)=\langle \l F+(1-\l)G , \mu\rangle - \int \int _{\R^{2d}}\log \E_{x,y}\left[e^{\l (F-\tilde W)+(1-\l)(G-\tilde W)}\right]\\
  &>\l\left(\langle  F, \mu\rangle - \int \int _{\R^{2d}}\log \E_{x,y}\left[e^{F-\tilde W }\right]\right)+(1-\l)\left(\langle G, \mu\rangle - \int \int _{\R^{2d}}\log \E_{x,y}\left[e^{G-\tilde W }\right]\right).
  \end{aligned}  
\end{equation}
 Furthermore, we see that $F\equiv 0$ solves \eqref{variationaltrap}. The variational equations yield for any $h \in \Ccal_\beta$.
\begin{equation}
    \langle h, \mu \rangle = \int_{\R ^d}\int_{\R ^d}\frac{\E_{x,y} [he^{-\Tilde{W}}]}{\E_{x,y} [e^{-\Tilde{W}}]}q(\d x , \d y),
\end{equation}
which identifies $\mu^*_q$.
\end{proofsect}
\qed

\begin{cor}\label{corproofmin}
    There exists a measure $q^*\in \Mcal_1^{(s)}(\R ^d\times\R ^d)$ such that 
    \begin{equation}\label{minimizerpairmeasure}
        \begin{aligned}
            &H(q^*|\overline{q^*}\otimes \m) -\int_{\R^ d} \int_{\R ^d}\log \E_{x,y}[e^{-{\tilde W }}]q^*(\d x,\d y)-\langle \log g,q^*\rangle \big \}\\
            &=\inf_{q\in \Mcal^{(s)}_1(\R ^ d\times \R ^d)}\left\{H(q|\overline{q}\otimes \m) -\int_{\R^ d} \int_{\R ^d}\log \E_{x,y}[e^{-{\tilde W }}]q(\d x,\d y)-\langle \log g,q\rangle \right \}.
        \end{aligned}
    \end{equation}
\end{cor}
\begin{proofsect}{Proof}
Recall that the mapping $q \mapsto H(q|\overline{q}\otimes\m$) has compact level sets, and thus, if $q_N$ is a minimizing sequence of the right hand side of \eqref{infIg}, then there exists a $q^*\in \Mcal_1^{(s)}(\R^d\times \R ^d)$ such that $q_N\Rightarrow q^*$ and $q^*$ minimizes the right hand side of \eqref{infIg}, as the mappings $q \mapsto \langle\log g ,q\rangle $ and $q \mapsto \langle\log \E_{x,y}\big[e^{-\Tilde{W}}\big] ,q\rangle$ are lower semi-continuous due to the fact that the mappings $(x,y)\mapsto \log p_\beta(x,y)$ and $(x,y) \mapsto \log \E_{x,y}\big[e^{-\Tilde{W}}\big]$ are continuous in $\R^ d\times \R ^d$.

\qed
\end{proofsect}

\subsection{Proof of Corollary~\ref{eigencor}}
     To prove (i), we first show that $$\inf\limits_{q\in\Mcal_1^{\ssup{{\rm s}}}(\R^d\times\R^d)}\Big\{H(q|\overline{q}\otimes\m)-\int_{\R ^d}\int_{\R ^d}\log \E^\beta_{x,y}\Big[e^{-\Tilde{W}}\Big]q(\d x, \d y)-\langle g,\log g\rangle\Bigr\}\leq -\log \lambda_T.$$
     Considering on the right hand side of \eqref{infIg} $$q(\d x , \d y)=\frac{1}{\lambda_T}\varphi_T(x)\varphi_T(y)\E_{x,y}^{\beta}\left [e^{-\Tilde{W} }\right]g(x,y)\m(\d y)\m( \d x).$$
    Then it is easy to see that $q$ is symmetric and $\overline{q}(\d x)=\varphi^2_T(x)\m(\d x)$, thus 
    $$\frac{q(\d x, \d y)}{\overline{q}(\d x)\m(\d y)}=\frac{1}{\lambda_T}\frac{\varphi_T(y)}{\varphi_T(x)}\E_{x,y}^{\beta}\left [e^{-\Tilde{W}}\right]g(x,y).$$
    Since $q$ has equal marginals, we have that 
    \begin{equation}\label{marginalproperty1}
     \int_{\R ^d }\int_{\R ^d }\log\frac{\varphi_T(y)}{\varphi_T(x)}q (\d x, \d y)=\int_{\R ^ d }\log\varphi_T(y) \overline{q}(\d y)-\int_{\R ^ d }\log\varphi_T(x) \overline{q}(\d x)=0,   
    \end{equation}
    and therefore the inequality follows. \smallskip

    \noindent To prove the converse inequality we write the infimum of the right hand side of \eqref{infIg} as
    \begin{equation}\label{equation1}
    \inf_{q \in \Mcal_1^{(s)}(\R ^d \times \R ^d)}\left \{\int_{\R ^d }\int_{\R ^d }\log\frac{q(\d x , \d y)}{\overline{q}(\d x)\E_{x,y}^{\beta}\left [e^{-\Tilde{W} }\right]g(x,y)\m(\d y)}q(\d x, \d y)\right\}.    
    \end{equation}
    Using the marginal property of $q$ as seen in \eqref{marginalproperty1}, we write \eqref{equation1} as 
    \begin{equation}
        \inf_{q \in \Mcal_1^{(s)}(\R ^d \times \R ^d)}\left \{\int_{\R ^d }\int_{\R ^d }\log\frac{q(\d x , \d y)}{\overline{q}(\d x)\frac{\varphi_T(y)}{\varphi_T(x)\lambda_T}\E_{x,y}^{\beta}\left [e^{-\Tilde{W} }\right]g(x,y)\m(\d y)}q(\d x, \d y)-\log\lambda_T\right\}.
    \end{equation}
    Note that $\frac{\varphi_T(y)}{\varphi_T(x)\lambda_T}\E_{x,y}^{\beta}\left [e^{-\Tilde{W} }\right]g(x,y)\m(\d y)$ is a probability measure for all $x \in \R ^d$ therefore the quantity 
    $$\int_{\R ^d }\int_{\R ^d }\log\frac{q(\d x , \d y)}{\overline{q}(\d x)\frac{\varphi_T(y)}{\varphi_T(x)\lambda_T}\E_{x,y}^{\beta}\left [e^{-\Tilde{W} }\right]g(x,y)\m(\d y)}q(\d x, \d y)$$
    is a relative entropy, and thus non-negative, which leads to the desired result. 

    \noindent (ii) follows directly from (ii) of Theorem \ref{unmin}.
    \qed

\subsection{Proofs of Theorems~\ref{propmininmizer} and \ref{mintheorem}.} Let us make an overview of the logic of our proofs. In the definition of the measure $\mu^*_{g,\m,W}$ in \eqref{minimizergen}, the pair measure $q^* \in\Mcal_1^{\ssup{{\rm s}}}(\R^d\times\R^d)$ depends on a measure $\m \in \Mcal(\R ^d)$ and on a function $g:\R ^ d \times \R ^d \to (0,\infty)$. Since $\mu^*_{g,\m,W}$ depends on $q^*$, as well as on a path density of the form $e^{-\tilde W(\omega)}$ (with respect to a path measure as seen (ii) Theorem~\ref{unmin}), it is reasonable to study some class of functions $W$ defined on $\Ccal_\beta$ and pair measures $q \in\Mcal_1^{\ssup{{\rm s}}}(\R^d\times\R^d)$ such that the probability measure of the form
\begin{equation}\label{genprobschr}
    \int_{\R ^ d}\int_{\R ^ d} \frac{e^{-\tilde W (\omega)}\P^{\beta}_{x,y}(\d \omega)}{\E^{\beta}_{x,y}[e^{\tilde W }]}q(\d x , \d y), \qquad \omega \in \Ccal_\beta,
\end{equation}
is a Schr\"odinger process, i.e., it solves the optimal transport problem $$ \inf\{H(Q|\P_\m): Q\circ\pi_0^{-1}=\nu_1,Q\circ\pi_\beta^{-1}=\nu_2\},$$ for some fixed probability measures $\nu_1,\nu_2$ on $\R ^d$, where $\P_{\m}$ is the Wiener measure with initial distribution $\m$. According to \cite{FG97} Q is a Schr\"odinger process if and only if $\frac{\d Q}{\d \P_{\m}}(\omega)=S_1(\omega(0))S_2(\omega(\beta))$, where $S_1$ and $S_2$ belong to $L^1(\R ^d)$. In other words $Q=\int_{\R ^ d}\int_{\R ^ d} \P^\beta_{x,y}S_1(x)S_2(y)\nu (\d x, \d y)$, where $\nu (\d x, \d y)$ is the joint distribution of the initial and final time of the Brownian Motion. In our case however, we initially see that the density of $\mu^*_{g,\m,W}$ depends on the entire path rather than the initial and terminal site. As we will see, if we properly choose the initial distribution $\m$, the minimizing pair measure $q^*$ and the function $g$, then the path dependence vanishes and the transformed Brownian motion in \eqref{genprobschr} is indeed a Schr\"odinger process. Let us briefly explore these conditions. Denote by $\l(g,\m,W)$ the quantity 
\begin{equation}\label{qminF}
\inf\limits_{q\in\Mcal_1^{\ssup{{\rm s}}}(\R^d\times\R^d)}\Big\{H(q|\overline{q}\otimes\m)-\int_{\R ^d}\int_{\R ^d}\log \E^\beta_{x,y}\Big[e^{\tilde W  }\Big]q(\d x, \d y)-\langle q,\log g\rangle\Bigr\}, 
\end{equation}
when the latter is finite. If the assumptions of Corollary~\ref{eigencor} are satisfied, then there exists a $\m$- square normalized function $\varphi:\R^d\to [0,\infty)$ such that for all $x \in \R ^d$, 
\begin{equation}\label{condpsi}
    e^{-\l(g,\m,W)}\int_{\R^d}\varphi(y)\E^\beta_{x,y}\Big[e^{\tilde W }\Big]g(x,y)\m(\d y)=\varphi(x).
\end{equation}
\noindent From the proof of Corollary~\ref{eigencor}, the minimizer of \eqref{qminF} is given by 
\begin{equation}\label{q*identtrapgen}
    q^*(\d x , \d y)=e^{-\l(g,\m,W)}\varphi(x)\varphi(x)\E^\beta_{x,y}\Big[e^{\tilde W  }\Big]g(x,y)\m(\d x)\m(\d y),
\end{equation}
and the measure in \ref{genprobschr} reads
\begin{equation}\label{mu*identF}
    e^{-\l(g,\m,W)}\int_{\R ^ d}\int_{\R ^ d} e^{\tilde W(\omega)}\P^{\beta}_{x,y}(\d \omega)g(x,y)\varphi(x)\varphi(x)\m(\d x)\m(\d y).
\end{equation}
For a Brownian motion, with initial measure $\m$, the joint distribution of the marginals at the initial time and final times are given by the measure $p_{\beta}(x,y)\m(\d x) dy$. We therefore require that the initial measure $\mu$ be absolutely continuous and for $g$ be a multiple of $p_{\beta}$. We also explore two cases for the trap potential $W$: (i) the Hard-wall potentials (path functions of the form $F(\omega)=0$ if $\omega[0,\beta]\subset \L$ and $F(\omega)=-\infty$, otherwise, where $\L\subset \R ^d$ is a finite box). Note that if the trap potential $W$ is of the form $W=\infty\1_{\L^c}$, then $F(\omega):=-\frac{1}{\beta}\int_0^\beta W(\omega(s))\d s$ is a hard wall path function, and (ii) the Soft-wall potentials which are continuous, finite everywhere path functions, of the form $\int_0^\beta W(\omega_s)ds$ for a suitably behaved $W:\R^d \to [0,\infty)$ (see Theorem~\ref{mintheorem}). for the case (i), we pick $\m=\Leb_\L$ and $g\equiv p_{\beta},$ and in (ii), we pick $\m=\Leb$ and again $g\equiv p_\beta$. \bigskip
\subsection{Proof of Theorem~\ref{propmininmizer}.} Here, we set the initial measure $\m$ to be the Lebesgue measure on $\L$, where $\L$ is a compact and connected subset of $\R^d$, and $g\equiv p_{\beta}$. 

\noindent It is easy to see that the minimizer in \eqref{genprobschr} is of the form  
\begin{equation}\label{mu^*identhw}
    \mu_{q_{\L}}^*=\int_{\L}\int_{\L}\P_{x,y}^{\beta,\L}q_{\L}(\d x,\d y),
\end{equation}
where $\P_{x,y}^{\beta,\L}$ is the normalized Brownian Bridge measure starting from $x\in \L$ and terminating at $y\in \L$ and remains in $\L$ for all $t\in [0,\beta]$. Additionally, $q_{\L}$ is the minimizer of \eqref{qminF} when $\m=\Leb_{\L}$, $g \equiv p_\beta$, i.e., the minimizer of
\begin{equation}\label{qminFhw}
    \inf\limits_{q\in\Mcal_1^{\ssup{{\rm s}}}(\R^d\times\R^d)}\Big\{H(q|\overline{q}\otimes\Leb_{\L})-\langle q,\log p_{\beta,\L}\rangle\Bigr\},
\end{equation}
where for $x,y\in \L$, $ p_{\beta,\L}(x,y):=\mu^{\beta}_{x,y}(\Ccal_\beta;B[0,\beta]\subset\L)$ is the Gaussian density in $\L$ (see Section 2.1 in \cite{CZ95} for the explicit formula of $p_{\beta,\L}$ and its properties). It is easy to see that $q_{\L}(\d x,\d y):=\frac{1}{|\L|}p_{\beta,\L}(x,y)\d x \d y,$ $x,y\in \L$ minimizes \eqref{qminFhw} as $q_\L$ is symmetric and $\overline{q_\L}(\d x)=\frac{1}{|\L|}\d x$, $x\in \L$. The measure \eqref{mu^*identhw} becomes $\mu_{q_\L}^*=\frac{1}{|\L|}\int_{\L}\int_{\L}\P_{x,y}^{\beta,\L}p_{\beta,\L}(x,y)\d x \d y,$ $x,y\in \L$, which is the desired form for a Schr\"odinger process. Moreover, note that the latter is a solution to the SDE 
\begin{equation}\label{sdehw}
    \d  X_t=\d B_t^{\L}, \qquad X_0\text{ is uniformly distributed on } \L,
\end{equation}
where $B_t^\L$ is the Brownian motion on $\L$. In particular, $\mu^*$ is ergodic with invariant measure being the uniform distribution on $\L$.
\subsection{Soft-wall potentials- proof of Theorem \ref{mintheorem}.}\label{sec-sch}Here we drop the subscripts $\m$ and $g$ from the rate function $I_{\m,W,g}^{(sym)}$. Let $\mu \in \Mcal_1(\Ccal_\beta)$ be such that
\begin{equation}\label{rtfnctF}
I_W^{(sym)}(\mu):=\inf_{q\in \Mcal^{(s)}_1(\R ^ d\times \R ^d)}\big\{H(q|\overline{q}\otimes \d x)+I_W^\ssup{q}(\mu)-\langle \log p_\beta,\mu _{0,\beta}\rangle \big\}
\end{equation}
is finite. 
According to Proposition~\ref{genprobschr} we have
\begin{equation}\label{minSWF}
        \inf_{\mu \in \Mcal_1(\Ccal_\beta)}\{I_{W}^\ssup{sym}(\mu)\}=\inf_{q\in \Mcal^{(s)}_1(\R ^ d\times \R ^d)}\Big\{H(q|\overline{q}\otimes \d x) -\int_{\R^ d} \int_{\R ^d}\log \E_{x,y}[e^{-\tilde{W}}]q(\d x,\d y)-\langle \log p_\beta,q\rangle \Big \},
\end{equation}
 and that the unique minimizer of \eqref{minSWF} is of the form  
 \begin{equation}\label{genprobschrF}
    \mu^*_{W}(\d \omega)=\int_{\R ^ d}\int_{\R ^ d} \frac{e^{-\tilde{W}(\omega)}\P^{\beta}_{x,y}(\d \omega)}{\E^{\beta}_{x,y}[e^{-\tilde{W} }]}q^*(\d x , \d y), \qquad \omega \in \Ccal_\beta,\quad q^*\in \Mcal^{(s)}_1(\R ^ d\times \R ^d),
\end{equation}
  where $q^*$ minimizes the right hand side of \eqref{minSWF}. From the assumptions regarding Soft-Wall potentials, there exists a normalized positive $\phi \in H^2(\R^d)$ solution to the problem $\Delta \phi +W\phi=\l(W)\phi$, where 
\begin{equation}\label{eigenvdefiF}
\l(W)=\sup_{\varphi\in H^1(\R^d)\colon , \|\varphi\|_2=1}\Big(\langle W,\varphi^2\rangle-\|\nabla\varphi\|_2^2\Big). 
\end{equation}
We also consider the following diffusion
 \begin{equation}\label{sdeshF}
     \d X_t=\frac{\nabla \varphi}{\varphi}\d t+\d B_t,\quad X_0 \text{ has distibution with density } \varphi^2 \text{ with respect to the Lebesgue measure}.
 \end{equation}
Then, the solution $\mu_W$ to the \eqref{sdeshF} is ergodic with invariant measure $\phi^2(x)\d x$. According to Girsanov's Theorem (see \cite[Prop. VIII.3.1]{RY99}), we have that $\d \mu_W =D_t^{W}\d \P_{\phi}$, where 
\begin{equation}\label{Girztrans}
    D_\beta^{(W)}:=e^{\int_0^\beta W(B_t)\d t-\beta \l(f)}\frac{\phi(B_\beta)}{\phi(B_0)},
\end{equation}
and $\P_\phi$ is the Wiener measure, starting with distribution $\phi^2(x)\d x$.

\noindent We now show that if $q^*(\d x, \d y)=\phi (x)\phi(y)\E_{x,y}\big[e^{-\int_0^\beta W(B_t)\d t-\beta \l(f)}\big]p_\beta(x,y)\d x\, \d y$, then the measure $\mu^*_W$ in \eqref{genprobschrF} is a Markov process. To see this, first we show that $q^*$ is a probability measure on $\R ^d\times \R ^d$ with equal marginals. Indeed, note that $q^*$ is symmetric and that its first marginal is given by
 \begin{equation}\label{marginalq}
     \overline{q^*}(\d x)=\phi^2(x)\Big(\int_{\R^d}\frac{\E_x\big[D_\beta^{(W)};B_\beta \in \d y\big]}{\d y}\d y\Big)\d x=\phi^2(x)\E_x\big[D_\beta^{(W)}\big]\d x.
 \end{equation}
 Note that $(D_t^{(W)})_{t\geq 0}$ is a Martingale under $\P_x$ (The Wiener measure starting from $x\in \R ^d$) with respect to the canonical Brownian Filtration. Thus $\overline{q^*}(\d x)=\phi^2(x)\d x$. Thus, $q^*$ is a probability measure on $\R ^d \times \R ^d$ with equal marginals. Now let $A\subset \Ccal_\beta$ be a Borel subset. Using (A.1) in the Appendix in \cite{S98} We have that 
 \begin{equation}
    \begin{aligned}
        &\mu^*_W(A)=\int_{\R ^d}\int_{\R^ d} \phi(x)\phi(y)\E_x[\1_Ae^{-\tilde{W}-\beta\lambda_\Lambda(W)}\delta_y(B_\beta)\d y]\d x=\\
        &=\int_{\R^d}\E_x[\1_AD_\beta^{(W)}]\phi^2(x)\d x=\mu_W(A).
        \end{aligned}
    \end{equation}
 Thus, $\mu^*_W=\mu_W$. 

\noindent We now show that $q^*(\d x, \d y)=\phi (x)\phi(y)\E_{x,y}\big[e^{\int_0^\beta W(B_t)\d t-\beta \l(f)}\big]p_\beta(x,y)\d x\, \d y$ minimizes the right hand side of \eqref{minSWF}. Note that we can rewrite the right hand side of \eqref{minSWF} as
 \begin{equation}\label{qminleb}
    \inf_{q\in \Mcal^{(s)}_1(\R^d\times \R ^d )}\Big\{\int_{\R^d}\int_{\R^d}
    \log \frac{\d q(x,y)}{\d \overline{q}(x) \E_x[e^{-\tilde{W}
    -\beta\lambda(f)}\frac{\phi(B_\beta)}{\phi(B_0)};B_\beta \in dy]}\d q (x,y)\Big\}-\beta\lambda(f),
\end{equation}
since $\int_{\R^d}\int_{\R^d}\log \frac{\phi(x)}{\phi(y)}q(\d x,\d y)=0$ because of the property of equal marginals of $q$. Note also that $ \E_x[e^{-\tilde{W}
-\beta\lambda(f)}\frac{\phi(B_\beta)}{\phi(B_0)};B_\beta \in dy]$ is a probability measure for any $x\in \R ^d$. Thus, 
    $$\int_{\R^d}\int_{\R^d}
    \log \frac{\d q(x,y)}{\d \overline{q}(x) \E_x[e^{-\tilde{W}
    -\beta\lambda(f)}\frac{\phi(B_\beta)}{\phi(B_0)};B_\beta \in dy]}\d q (x,y)$$
is a relative entropy for any $q\in \Mcal^{(s)}_1(\R^d\times \R ^d )$, and thus it is non negative. Thus the expression in \eqref{qminleb} is not less than $-\beta \l(f)$. Now, if we set $q=q^*$, it is easy to see that $$\int_{\R^d}\int_{\R^d}
    \log \frac{\d q^*(x,y)}{\d \overline{q^*}(x) \E_x[e^{-\tilde{W}
    -\beta\lambda(f)}\frac{\phi(B_\beta)}{\phi(B_0)};B_\beta \in dy]}\d q^* (x,y)=0,$$
 Thus the right hand side of \eqref{qminleb} is not greater than $-\beta\l(f)$. Thus, $q^*$ minimizes the right hand side of \eqref{minSWF}, and Theorem \ref{mintheorem} follows.
 \qed
\subsection{Identification of the rate function $J_{\Leb,N,W}$- Proof of Theorem \ref{thDVtrap}}\label{sec-DV}\bigskip

\noindent Recall the rate function given by \eqref{Jnonnormal} when $\m=\Leb$ and when $g\equiv p_\beta$. We write 
\begin{equation}\label{ratelebgaussian}
      \overline{J}(p)=\inf\limits_{q\in\Mcal_1^{\ssup{{\rm s}}}(\R^d\times\R^d)}\Big\{H(q|\overline{q}\otimes\d x)+\Tilde{J}_W^{\ssup{q}}(p)\Big\},
\end{equation}
where
\begin{equation}\label{FLtranstrapdxgaussian}
    \Tilde{J}_W^{\ssup{q}}(p):=\sup_{f\in\Ccal_{\rm b}(\R^d)}\Bigl\{\beta\langle f,p\rangle -\int_{\R^d}\int_{\R^d} q(\d x,\d y)\,\log\frac{\E_{x}\Bigl[{\rm e}^{\int_0^{\beta}f(B_s)\,\d s};B_\beta \in \d y\Bigr]}{\d y}+\beta \langle W, p\rangle\Bigr\}.
\end{equation}

\noindent We split the proof of Theorem \ref{thDVtrap} in 2 steps. The first step is the proof of the lower bound $\overline{J}(p)\ge \beta I_W(p)$, where $I_W(p)$ is defined in \eqref{DVtrap}. In the second step, we show that 
\begin{equation}\label{upperfreebm}
    J^{\ssup{\rm sym}}_{\Leb}(p):=\inf\limits_{q\in\Mcal_1^{\ssup{{\rm s}}}(\R^d\times\R^d)}\Big\{H(q|\overline{q}\otimes\d x)+J^{\ssup{q}}(p)-\langle q,\log p_\beta\rangle\Big\}\leq \|\nabla \varphi\|_2^2,
\end{equation}
where $\frac{\d p}{\d x}:=\phi^2(x)$ and $\phi \in H^1(\R^d)$. This, implies the upper bound $\overline{J}(p)\leq \beta I_W(p)$.

\noindent Let us turn into details. for $f\in C_b(\R^d)$, let $\varphi_f$ be the unique positive $L^2$ normalized eigenfunction of $-\Delta+f-W$ with corresponding eigenvalue $\l(f-W)$ given by \eqref{eigenvdefi}. Then,
\begin{equation}\label{MGtrapplusf}
    D^{(f)}_\beta:=e^{\int_0^\beta (f-W)(B_s)\d s+\beta\,\l(f-W)}\frac{\varphi_f(B_\beta)}{\varphi_f(B_0)} 
\end{equation}
defines a martingale $(D^{(f)}_\beta)_{\beta \ge 0}$ under $\P_x$ for any $x \in \R ^d$ with respect to the canonical Brownian filtration (see Section \ref{sec-sch} for details). Our aim is to construct $D^{(f)}_\beta$ within $\Tilde{J}_W^{\ssup{q}}(p)$. To do this we need the following Lemma
\begin{lemma}\label{L2lemma}
    If $\langle W,p \rangle<\infty$, then
    \begin{equation}\label{L2lemmaeq}
        \Tilde{J}_W^{\ssup{q}}(p)\geq \sup_{f\in\Ccal_{\rm b}(\R^d)}\Bigl\{\beta\langle f,p\rangle -\int_{\R^d}\int_{\R^d} q(\d x,\d y)\,\log\frac{\E_{x}\Bigl[{\rm e}^{\int_0^{\beta}(f-W)(B_s)\,\d s};B_\beta \in \d y\Bigr]}{\d y}\Bigr\}.
    \end{equation}
\end{lemma}
\begin{proofsect}{Proof}
    Let $f\in C_b(\R ^d)$ be given. For $M>0$, consider $f_M=(f+W)\land M$, then $f_M\in C_b(\R ^d)$ with $f_M\uparrow f+W$. Since $\langle W,p \rangle<\infty$, we have
    $$\limsup_{M\to \infty}\langle f_M,p\rangle \le \langle f+W,p\rangle. $$
    Moreover, by the Monotone Convergence Theorem we have that 
    $$\lim_{M\to \infty}\frac{\E_{x}\Bigl[{\rm e}^{\int_0^{\beta}(f_M-W)(B_s)\,\d s};B_\beta \in \d y\Bigr]}{\d y}=\frac{\E_{x}\Bigl[{\rm e}^{\int_0^{\beta}f(B_s)\,\d s};B_\beta \in \d y\Bigr]}{\d y},$$
    which finishes the proof of Lemma. \ref{L2lemma} 
    \qed
\end{proofsect}

\noindent Using Lemma \ref{L2lemma}, the marginal property of $q$ (i.e., $\int\int \log \frac{\varphi(y)}{\varphi(x)}q(\d x,\d y)=0$) and the fact that $-\langle W,p\rangle\leq0$, we have that 
\begin{equation}\label{GirLB1}
    \Tilde{J}_W^{\ssup{q}}(p)\ge \sup_{f\in\Ccal_{\rm b}(\R^d)}\Bigl\{\beta\big(\,\langle f-W,p\rangle+\l(f-W)\,\big) -\int_{\R^d}\int_{\R^d} q(\d x,\d y)\,\log\frac{\E_{x}\Bigl[D^{(f)}_\beta;B_\beta \in \d y\Bigr]}{\d y}\Bigr\}.
\end{equation}
We now substitute \eqref{GirLB1} to \eqref{ratelebgaussian} and we have that 
\begin{equation}\label{GirLB2}
    \overline{J}(p)\geq \inf\limits_{q\in\Mcal_1^{\ssup{{\rm s}}}(\R^d\times\R^d)}\sup_{f\in\Ccal_{\rm b}(\R^d)}\Big\{\beta\big(\,\langle f-W,p\rangle+\l(f-W)\,\big)+\int_{\R^d}\int_{\R^d}\log\frac{q(\d x,\d y)}{\overline{q}(\d x)\E_{x}\Bigl[D^{(f)}_\beta;B_\beta \in \d y\Bigr]}\Big\}.
\end{equation}
By the Martingale property of $ (D_{\beta}^{\ssup{f}})_{\beta\ge 0} $, the measure $ \E_x[D_{\beta}^{\ssup{f}};B_{\beta}\in\d y] $ is a probability measure on $ \R^d $ for any $ x \in\R^d $. Hence, the double integral in \eqref{GirLB2} is an entropy between probability measures and therefore non-negative, by Jensen's inequality. Thus, we have that
\begin{equation}\label{GirLB3}
    \overline{J}(p)\geq \inf\limits_{q\in\Mcal_1^{\ssup{{\rm s}}}(\R^d\times\R^d)}\sup_{f\in\Ccal_{\rm b}(\R^d)}\Big\{\beta\big(\,\langle f-W,p\rangle+\l(f-W)\,\big)\Big\}.
\end{equation}
Note that the map $f\mapsto \l(f-W) $ is the Legendre-Fenchel transform of  $ I_W $, as is seen from the Rayleigh-Ritz principle in \eqref{eigenvdefi}. According to the Duality Lemma \cite[Lemma~4.5.8]{DZ98}, the r.h.s.\ of \eqref{GirLB3} is equal to $ \beta I_\L(p) $ since it is equal to the Legendre-Fenchel transform of $\l(\cdot\,-W)$. Hence step 1 is complete

\noindent We now show the second step. We fist show that \eqref{upperfreebm} holds for any absolutely continuous and compactly supported measure $p\in \Mcal_1(\R^d)$. Let $ \varphi=\sqrt{\frac{\d p}{\d x}} $ with compact support $\L\subset \R ^d$. Then, it is easy to see that $J^{\ssup{\rm sym}}_{\Leb}(p)\leq J^{\ssup{\rm sym}}_{\L}(p)$, where
\begin{equation}\label{JsymLambdadef}
J^{\ssup{\rm sym}}_{\L}(p)=\inf_{q\in\Mcal_1^{\ssup{\rm s}}(\R^d\times\R^d)}\Big\{H(q|\overline{q}\otimes\Leb_\L)+ J^{\ssup{q}}_{p_\beta}(p)\Big\},\qquad p\in\Mcal_1(\R^d),
\end{equation}
and
\begin{equation}\label{JqLambdadef}
J^{\ssup{q}}_{p_\beta}(p)=\sup_{f\in\Ccal_{\rm b}(\R^d)}\Bigl\{\beta\langle f,p\rangle-\int_{\R^d}\int_{\R^d} q(\d x,\d y)\,\log\E_x\big[ {\rm e}^{\int_0^{\beta}f(B_s)\,\d s};B_\beta\in\d y\big]\big/\d y\Bigr\},
\end{equation}
as the set of the pair measures $q\in \Mcal_1^{\ssup{\rm s}}(\R^d\times\R^d)$ having compact support on $\L\times \L$ is subset of all the pair measures $q \in \Mcal_1^{\ssup{\rm s}}(\R^d\times\R^d)$ (see Lemma 2.1 in \cite{AK08} for details). From the proof of steps 2 and 3 in Theorem 1.5 in \cite{AK08} we readily have that $J^{\ssup{\rm sym}}_{\L}(p)\leq \beta \|\nabla \varphi\|^2_2 $, thus, the upper bound in \eqref{upperfreebm} holds for any absolutely continuous and compactly supported measure $p\in \Mcal_1(\R^d)$. 

\noindent If now $ \varphi=\sqrt{\frac{\d p}{\d x}} \in H^1(\R ^d) $, we can approximate $\varphi$ in the $H^1$ norm, by an $L^2$ normalized sequence of $\varphi_N$ of functions that are smooth and compactly supported on a set $\L_N$, and thus $\|\nabla \varphi_N\|^2_2\to \|\nabla \varphi\|^2_2$, as $N\to \infty$. finally, we have that $J^{\ssup{\rm sym}}_{\Leb}(p)\leq J^{\ssup{\rm sym}}_{\L_N}(p)$, for all $N\in \N$ and all $p\in \Mcal_1(\R ^d)$ with compact support on $\L_N$. Furthermore, since $J^{\ssup{\rm sym}}_{\Leb}$ is lower semi-continuous we have that $J^{\ssup{\rm sym}}_{\Leb}(\varphi)\leq \liminf_{N\to \infty}J^{\ssup{\rm sym}}_{\Leb}(\varphi_N)$. Using the above, the upper bound \eqref{upperfreebm} follows, since $$J^{\ssup{\rm sym}}_{\Leb}(\phi)\leq\liminf_{N\to \infty}J^{\ssup{\rm sym}}_{\Leb}(\varphi_N)\leq\liminf_{N\to \infty}J^{\ssup{\rm sym}}_{\L_N}(\varphi_N)=\liminf_{N\to \infty} \beta \|\nabla \varphi_N\|^2_2=\beta \|\nabla \varphi\|^2_2$$. 
\qed

\section{Appendix}
\subsection{Proof of the auxiliary results}\label{sec-auxresults} \smallskip
\subsubsection{Proof of Proposition~\ref{main}}We split the proof in three lemmas. The first lemma provide an explicit formula for the rate function $I_{W,\m}$, whereas the second and third lemma identify $I_{W,\m}$ with the one stated in Theorem~\ref{main}.

\begin{lemma}
Fix $ \beta \in (0,\infty) $. Then as $ N\to\infty$, under the measure $ \P_{\m,N,W}^{\ssup{{\rm sym}}} $ the empirical path measures $ L_N $ satisfy a Large Deviations Principle on $ \Mcal_1(\Ccal_\beta) $ with speed $ N $ and rate function $I_{W,\m}:\Mcal_1(\Ccal_\beta)\to \R \cup \{\infty\}$ given by
\begin{equation}\label{JTtrapRF}
\hat{I}_{W,\m}(\mu)=\begin{cases}
    H(\mu_0|\m)+H(\mu_{0,\beta}|\mu_{0}\otimes\mu_{\beta})+\Scal_W(\mu),\text{ if } \mu_0=\mu_\beta\\
    \infty, \text{ otherwise},
\end{cases}
\end{equation}
where $\Scal_W(\mu)=\sup_{F \in C_b(\Ccal_\beta)}\left\{\langle F,\mu\rangle - \int_{\R ^d}\int_{\R ^d}\E^{\beta}_{x,y}\left[e^{F -\int_0^\beta W(B_s)ds}\right]\mu_{0,\beta}(\d x, \d y) \right\}$.
\end{lemma}
\begin{proofsect}{Proof}
\noindent Following [T07], we apply Theorem 2.1 in [G02]. To this end, it suffices to show that $L_N$ under the measure $\P_{\m,N,W}^{\ssup{{\rm sym}}}$, is a mixture of Large Deviations Systems in the sense of Definition 2.1. in [G02]. For completeness we follow the same notation as in [T07], and we show that all the requirements of Definition 2.1 in [G02] are satisfied. To see this, we have
\begin{itemize}
    \item $\Zcal:=\Mcal_1(\Ccal_\beta)$ is a polish space when it is equipped with the topology of weak convergence of probability measures 
\item $\Xcal_\infty:=\Mcal_1(\R ^d \times \R^d)$ equipped with topology of weak convergence of probability measures 
    \item Let $N\geq 1$. We denote by $$\Xcal_N:=\{\frac{1}{N}\sum_{i=1}^N\delta_{(s_i,a_i)}\in \Mcal_1(\R^d\times\R^d):s_i,a_i \in \R ^d, i=1,\dots,N \},$$
    and for every $q\in \Xcal_\infty,$ there exists $q_N \in \Xcal_N$ such that $q_N \to q$ as $N \to \infty$, in the weak sense.
    \item the map $\pi : \Zcal \to \Xcal_\infty$, with $\pi(\mu)=\mu_{0,\beta}$ is continuous and surjective, as seen on [B68], Theorem 2.7, since the projection mapping $\pi_{0,\beta}:\Ccal_\beta \to \R^d \times \R^ d$ is continuous.
    \item For every $N\ge 1$, and every $q_N=\frac{1}{N}\sum_{i=1}^N\delta_{(s_i,a_i)} \in \Xcal_N$, we denote by $P_{q_N}$ the distribution of $L_N$ under the measure $\bigotimes_{i=1}^Ne^{-\sum_{i=1}^N\int_0^\beta W(B^{(i)}_s)ds}\P^\beta_{(s_i,a_i)}$. Moreover, the family $$\Pi:=\{P^N_{q_N}:q_N \in \Xcal_N , N\geq 1\}$$ is such that $$P^N_{q_N}(\pi^{-1}(q_n)^c)=\int_{\{\pi(L_N)\neq \frac{1}{N}\sum_{i=1}^N\delta_{(s_i,a_i)}\}}\bigotimes_{i=1}^Ne^{-\sum_{i=1}^N\int_0^\beta W(\omega^{(i)}_sds}\P^\beta_{(s_i,a_i)}(\d \omega)$$ $$\leq \bigotimes_{i=1}^N\P^\beta_{(s_i,a_i)}\left(\pi(L_N)\neq \frac{1}{N}\sum_{i=1}^N\delta_{(s_i,a_i)}\right)=0 $$
    \item Denote by $Q_N$ the distribution of $\frac{1}{N}\sum_{i=1}^N\delta_{(B^{(i)}_0,B^{(\sigma_N(i))}_0)}$, where $\sigma_N$ is a sequence of random variables uniformly distributed on $\Sym_N$. Then for any $A \in \Zcal$, $A$ Borel set, we have
    $$\P_{\m,N,W}^{\ssup{{\rm sym}}}(L_N \in A)=\int_{\Xcal_N}P_{q}(A)dQ_N(q).$$
\end{itemize}
Moreover, according to [T07], Theorem 2, we have that the sequence $\frac{1}{N}\sum_{i=1}^N\delta_{(B^{(i)}_0,B^{(\sigma_N(i))}_0)}$ satisfies a LDP on $\Mcal_1(\R^d\times \R ^d)$ with speed $N$ and good rate function 
\begin{equation}\label{JTratefunction}
    \Ical(\nu)=\begin{cases}
        H(\overline{\nu}|\m)+H(\nu|\overline{\nu}\otimes\overline{\nu}), \text{ if } \int_{\R ^d}\nu(\d x, \cdot)=\int_{\R ^d}\nu(\cdot, \d y)=\overline{\nu}\\
        \infty, \text{ otherwise}
    \end{cases}
\end{equation}
Finally, let $ q \in \Xcal_\infty$ and let $q_N \in \Xcal_N$, such that $q_N=\frac{1}{N}\sum_{i=1}^N\delta_{(s_i,a_i)}$, and $q_N\to q$ in the weak sense. We show that $L_N$ satisfies a LDP with speed $N$ and good rate function $\Scal_W$. Indeed, let $F \in C_b(\Ccal_\beta)$ and let 
$$\Lambda_N(F):=\frac{1}{N}\log E_{q_N}[e^{-N\langle F,L_N\rangle}],$$
where $E_{q_N}$ is the expectation with respect to $P_{q_N}$. Then since for any $\omega \in \Ccal_\beta$ $\omega \mapsto e^{F(\omega)-\int_0^\beta W(\omega_s)ds} $ is continuous and bounded, we have 
$$\frac{1}{N}\log E_{q_N}[e^{-N\langle F,L_N\rangle}]=\langle \log\E^\beta_{x,y}\left[e^{F -\int_0^\beta W(B_s)ds}\right],q_N\rangle \to\langle \log\E^\beta_{x,y}\left[e^{F -\int_0^\beta W(B_s)ds}\right],q\rangle$$
$$=\int_{\R ^d}\int_{\R ^d}\log\E^\beta_{x,y}\left[e^{F -\int_0^\beta W(B_s)ds}\right]q(\d x, \d y),$$
as $N\to \infty$. 
It easy to see that the map $F \mapsto \int_{\R ^d}\int_{\R ^d}\log\E^\beta_{x,y}\left[e^{F -\int_0^\beta W(B_s)ds}\right]q(\d x, \d y)$ is Gateaux differentiable. Additionally, exponentially tightness follows from the fact that $W$ is non-negative and according to Theorem 3 in [T07], $L_N$ is exponentially tight under the measure $\bigotimes_{i=1}^N\P^\beta_{(s_i,a_i)}$. The final result follows from \cite{DZ98}, Corollary 4.6.14, and from [G02], Theorem 2.1.

\qed

\end{proofsect}

\noindent In the following Lemma, we provide an alternative expression for \eqref{JTratefunction}, similar to the expression of the rate function in Proposition~\ref{main}. This expression is useful for the minimization of the rate function.

\begin{lemma}\label{lemmaalternativetrap}
    Define 
    \begin{equation}\label{Itildesymdef}
\widetilde{I}^{\ssup{{\rm sym}}}_{\m,W}(\mu)=\inf_{q\in\Mcal_1^{\ssup{\rm s}}(\R^d\times\R^d)}\Bigl\{H(q|\overline{q}\otimes\m)+ \widetilde{I}_W^{\ssup{q}}(\mu)\Bigr\},\qquad \mu\in\Mcal_1(\Ccal_\beta),
\end{equation} 
where
\begin{equation}\label{Jtildeqdef}
\widetilde{I}_W^{\ssup{q}}(\mu)=\sup_{F \in\Ccal_{\rm b}(\Ccal)}\Bigl\{\langle F ,\mu\rangle -\int_{\R^d}\int_{\R^d} \log\E_{x,y}^{\beta}\big[{\rm e}^{F -\tilde{W} }\big]q(\d x, \d y)\Bigr\},\qquad\mu\in\Mcal_1(\Ccal_\beta).
\end{equation}
Then $\hat{I}_{W,\m}\equiv \widetilde{I}^{\ssup{{\rm sym}}}_{\m,W}$
\end{lemma}
\begin{proofsect}{Proof}
    To see this, we relax the supremum in \eqref{Jtildeqdef} to functions of the form $f(\omega(0),\omega(\beta)$. Then we have 
    \begin{equation}
    \begin{aligned}
        &\sup_{F \in\Ccal_{\rm b}(\Ccal_\beta)}\Bigl\{\langle F ,\mu\rangle -\int_{\R^d}\int_{\R^d} \log\E_{x,y}^{\beta}\big[{\rm e}^{F -\tilde{W} }\big]q(\d x, \d y)\Bigr\}\\
        &\ge \sup_{f \in\Ccal_{\rm b}(\R^d\times \R ^d)}\Bigl\{\langle f ,\mu_{0,\beta}-q\rangle -\int_{\R^d}\int_{\R^d} \log\E_{x,y}^{\beta}\big[{\rm e}^{-\tilde{W} }\big]q(\d x, \d y)\Bigr\} \ge \sup_{f \in\Ccal_{\rm b}(\R^d\times \R ^d)}\Bigl\{\langle f ,\mu_{0,\beta}-q\rangle\Bigr\}.
    \end{aligned}
        \end{equation}
        Additionally, noting that  $$H(\mu_{0,\beta}|\mu_0\otimes\m)=H(\mu_{0,\beta}|\mu_0\otimes\mu_\beta)+H(\mu_{0}|\m),$$
        we have the stated result.
        \qed

\end{proofsect}
\begin{lemma}
    Recall the definition of $I_{W,\m}$ in \eqref{Trapratefuncemp}. Then we have $I_{W,\m}\equiv \widetilde{I}^{\ssup{{\rm sym}}}_{\m,W}(\mu)$ 
\end{lemma}
\begin{proofsect}{Proof}
    Recall the definition of $I^{\ssup{q}}(\mu)$ in \eqref{Jqdef}. It suffices to show that $\widetilde{I}_W^{\ssup{q}}(\mu)=I^{\ssup{q}}(\mu)+\langle \tilde{W},\mu\rangle$ for every $q \in \Mcal_1^{(s)}(\R^d\times \R^d)$. We first show that $\widetilde{I}_W^{\ssup{q}}(\mu)\ge I^{\ssup{q}}(\mu)+\langle \tilde{W},\mu\rangle$. To show this, let $G\in C_\b(\Ccal_\beta)$, and consider in \eqref{Jtildeqdef} $F=G+\tilde{W}_M$, where for $M>0$ we denote by $\tilde{W}_M:=\tilde{W}\land M$. Then we have
    \begin{equation}
    \begin{aligned}
    &\widetilde{I}_W^{\ssup{q}}(\mu)\ge \langle G+W_M ,\mu\rangle -\int_{\R^d}\int_{\R^d} \log\E_{x,y}^{\beta}\big[{\rm e}^{G +\tilde{W}_M -\tilde{W} }\big]q(\d x, \d y)\\
   & \ge \langle G ,\mu\rangle -\int_{\R^d}\int_{\R^d} \log\E_{x,y}^{\beta}\big[{\rm e}^{G }\big]q(\d x, \d y) +\langle \tilde{W}_M,\mu\rangle.
    \end{aligned}
    \end{equation}
    Here, we used the fact that $W_M-W\leq 0$. Using the Monotone Convergence theorem, we pass to the limit as $M\to \infty$ and the inequality follows. To prove the converse inequality consider in \eqref{Jqdef} $F=G-\tilde{W}_M$. Then we have
    \begin{equation}
        \begin{aligned}
            I^{\ssup{q}}(\mu)+\langle \tilde{W},\mu\rangle\ge \langle G+(W-W_M)  ,\mu\rangle -\int_{\R^d}\int_{\R^d} \log\E_{x,y}^{\beta}\big[{\rm e}^{G -\tilde{W}_M }\big]q(\d x, \d y) \\
            \ge \langle G  ,\mu\rangle -\int_{\R^d}\int_{\R^d} \log\E_{x,y}^{\beta}\big[{\rm e}^{G -\tilde{W}_M }\big]q(\d x, \d y).
        \end{aligned}
    \end{equation}
    Using again the Monotone Convergence Theorem we have the desired result.

    \qed
    \end{proofsect}

\subsubsection{Proof of corollary~\ref{main2}} To show this, we use Proposition~\ref{main}, and the Contraction principle \cite[Th.~4.2.1]{DZ98}. We consider the mapping $\Mcal_1(\Ccal_\beta)\ni\mu\mapsto\Psi(\mu):=\frac{1}{\beta}\int_0^\beta \mu\circ\pi_t^{-1}\d t$, where we recall that $ \pi_s(\omega)=\omega(s) $ is the projection. It is easy to see that $\Psi$ is continuous in the topology of the weak convergence of probability measures and that $\Psi(L_N)=Y_N$. Using the above we have that the sequence $Y_N$ satisfies the Large Deviations Principle under the measure defined in \ref{Psymdef} with rate function
\[\overline{J}_{\m,W}(p)=\inf_{\mu \in \Mcal_1(\Ccal_\beta):p=\Psi(\mu)}\{I^{(sym)}_{\m}(\mu)+\langle \Tilde{W},\mu\rangle-\inf_{\mu \in \Mcal_1(\Ccal_\beta)}\{I^{(sym)}_{\m}(\mu)+\langle \Tilde{W},\mu\rangle\}\}.\]
To show that $\overline{J}_{\m,W}=J_{\m,W}$, we recall that in [AK08] it has been proved that 
\[\inf_{\mu \in \Mcal_1(\Ccal_\beta):p=\Psi(\mu)}I^{(sym)}_{\m}(\mu)=J^{(sym)}_{\m}(p).\]
Moreover, using a change of variable argument we have that $\langle \Tilde{W},\mu\rangle=\langle {W},p\rangle$ for any $\mu \in \Mcal_1(\Ccal_\beta)$ such that $\Psi(\mu)=p$.

\qed

\subsubsection{Proof of Proposition~\ref{cormain}.}Assume first that $\m$ is a finite measure. Then the result is a direct application of \cite[Prop~1.3]{AK08} in conjunction with Theorem \ref{main} for (i), and Corollary \ref{main2} for (ii). Let now $\m$ be a measure with density $\psi$. We show that there exists a bounded continuous function $f_\m \in L^1(\R ^d)$ such that $\int_{\R ^d} f_\m(x)\m(\d x)=1$. To see this, consider the partition of $\R ^d$, $(A_n)_n$ with $A_n=\{x \in \R ^d : n-1\leq \psi(x)\leq n\}$, then $\m(A_n)<\infty$ and we set 
$$g:=\sum_{k=1}^\infty \frac{1}{2^k}\frac{\1_{A_k}}{1+\m(A_k)}.$$
then $0\leq g \leq 1$. We can now set $f_\m=\frac{g}{\int_{\R ^d}g(x)\m(\d x)}$, however, since $f_\m$ is not continuous we overcome this by setting 
$$g:=\sum_{k=1}^\infty \frac{1}{2^k}\frac{\phi_{A_k}(x)}{1+\m(A_k)},$$
where $\phi:\R ^d \to \R$ is a smooth approximation of $\1_{A_n}$. Thus since $f_{\m}\m (\d x)$ is a probability measure, we may apply the previous case with $g'\equiv g f_\m$ and our result follows.
\qed

\subsection{Large Deviations}\label{sec-LDP}\bigskip\\
In this chapter we will see an overview of some of the tools we will be using. The main concept of Large deviations theory is the asymptotic computation of small probabilities on an
exponential scale. On this essay we will have to study complicated and decaying families of probability measures, so a proper understanding of Large deviations theory is essential. Proofs of the results can be found in [DZ98] \\\\
Recall that a family of measures $(\mu_\epsilon)_{(\epsilon\geq 0)}$ on $(\mathcal{X},{\mathcal{B}}_\mathcal{X})$, where $\mathcal{X}$ is a topological space and ${\mathcal{B}}_\mathcal{X}$ is the Borel sigma-algebra on $\mathcal{X}$, satisfies the Large Deviation Principle (LDP) with speed $\frac{1}{\epsilon}$ and rate function $I:\mathcal{X}\rightarrow [0,\infty]$ if the level sets $L(a)=\{x \in \mathcal{X}: I(x)\leq a\}$ are closed in $\mathcal{X}$ for every $a \geq 0$ and the following holds:\\\\
Upper bound: For every $F\in \mathcal{B}_\mathcal{X}$ closed
\[\limsup_{\epsilon\rightarrow 0} \epsilon\log\mu_\epsilon(F) \leq -\inf_{x\in F}I(x) .\]
Lower bound: For every $G\in \mathcal{B}_\mathcal{X}$ open
\[-\inf_{x\in G}I(x)\leq \liminf_{\epsilon\rightarrow 0} \epsilon\log\mu_\epsilon(G).\]
Moreover, if the level sets are compact, we call $I$ good rate function.
\begin{example}
(Cramer's theorem) Let $\mathcal{X}=\R ^d$ and $(X_n)_{n\in \N}$ i.i.d $\R ^d$ valued random variables. Then the family defined as
\[\mu_n(A)=\P\Big(\frac{X_1+\dots+X_n}{n}\in A\Big)\quad A\in \mathcal{B}_{\R ^d},\]
satisfies the LDP with speed $n$ and rate function the Fenchel-Legendre transformation of the function $\lambda\mapsto \log \E [e^{\langle\lambda, X_1\rangle}]$, $\lambda \in \R ^d$ i.e.,
\[I(x)=\sup_{\lambda \in \R ^d}\{\langle \lambda ,x\rangle -\log \E [e^{\langle\lambda, X_1\rangle}]\}.\]
\end{example}

\begin{theorem}
(Contraction principle) Let $\Xcal$ and $\Ycal$ be Hausdorff topological spaces and $f:\Xcal \to \Ycal$ continuous function. If a  family of probability measures $\mu_\epsilon$ satisfies the LDP on $\Xcal$ with good rate function $I$ then the family $\nu_\epsilon= \mu_\epsilon \circ f^{-1}$ satisfies the LDP on $\Ycal$ with good rate function
\[J(y)=\inf_{x\in \Xcal: y=f(x)}\{I(x)\}.\]
\end{theorem}
\noindent We can also prove the inverse of the previous theorem under certain conditions, recall that a family of probability measures $\mu_\epsilon$ on $\Xcal$ is called exponentially tight if for every $a>0$ exists $K_a \subset \Xcal$ compact such that
\[\limsup_{\epsilon\rightarrow 0} \epsilon\log\mu_\epsilon(K_\epsilon) < -a\]
\begin{theorem}
(Inverse Contraction principle) Let $\Xcal$ and $\Ycal$ be Hausdorff topological spaces and $f:\Xcal \to \Ycal$ continuous and bijection function. If a  family of probability measures $\mu_\epsilon$ on $\Xcal$ exponentially tight and the family $\nu_\epsilon= \mu_\epsilon \circ f^{-1}$ satisfies the LDP on $\Ycal$ with good rate function $J$, then the family $\mu_\epsilon$ satisfies the LDP with good rate function $I=J\circ f$.
\end{theorem}
\begin{theorem}\label{vthm}
(Varadhan) Let $\Xcal$ be a Hausdorff regular space and Suppose that the family $\mu_\epsilon$ satisfies the LDP with a good rate function $I:\Xcal \to [0,\infty]$. Moreover, Let $Z_\epsilon$ be a family of $\Xcal$ valued random variables and let $\phi:\Xcal \to \R$ be any continuous function. Assume further either the tail condition
\begin{equation}\label{tail}
    \lim_{M\to \infty} \limsup_{\epsilon \to 0}\epsilon\log \E\big[e^{\frac{\phi(Z_\epsilon)}{\epsilon}}|\phi(Z_\epsilon)>M\big]=-\infty
\end{equation}
or the following moment condition for some $\gamma >1$
\begin{equation}\label{moment}
    \limsup_{\epsilon \to 0}\epsilon\log \E[e^{\frac{\gamma\phi( Z_\epsilon)}{\epsilon}}]<\infty
\end{equation}
then
\begin{equation}
    \lim_{\epsilon \to 0}\epsilon\log \E[e^{\frac{\phi( Z_\epsilon)}{\epsilon}}]=\sup_{x\in \Xcal}\{\phi(x)-I(x)\}
\end{equation}
\end{theorem}
\begin{remark}
The moment condition \ref{moment} implies the tail condition
\ref{tail}. Additionally, If $\phi$ is bounded from above, then tail condition
\ref{tail} is implied  
\end{remark}
\noindent Theorem \ref{vthm} is a direct consequence of the following two lemmas.
\begin{lemma}\label{limlo}
If $\phi:\Xcal \to \R$ is lower semi continuous and the Large Deviations lower bound holds with $I:\Xcal \to [0,\infty]$ then
\begin{equation}
    \liminf_{\epsilon \to 0}\epsilon\log \E[e^{\frac{\phi( Z_\epsilon)}{\epsilon}}]\geq\sup_{x\in \Xcal}\{\phi(x)-I(x)\}
\end{equation}
\end{lemma}
\begin{lemma}\label{lemup}
If $\phi:\Xcal \to \R$ is an upper semi continuous function for which
the tail condition \ref{tail} holds, and the Large Deviations upper bound holds with the good rate function $I:\Xcal \to [0,\infty]$, then
\begin{equation}
    \limsup_{\epsilon \to 0}\epsilon\log \E[e^{\frac{\phi( Z_\epsilon)}{\epsilon}}]\leq\sup_{x\in \Xcal}\{\phi(x)-I(x)\}
\end{equation}
\end{lemma}
\noindent We usually study measures of the form
\begin{equation}\label{normalized}
    \nu_\epsilon(A)=\frac{\int_Ae^{\phi(x)/\epsilon} \mu_\epsilon(\d x)}{Z_\epsilon(\phi)},\quad A\in \mathcal{B}_\Xcal
\end{equation}
where
\[Z_\epsilon(\phi)=\int_\Xcal e^{\phi(x)/\epsilon} \mu_\epsilon(\d x)\]
and $\mu_e$ is a family of probability measures on $\Xcal$ and $\phi:\Xcal \to \R$. We have the following corollary 
\begin{cor}
Let $\Xcal$ be a Hausdorff regular space and Suppose that the family $\mu_\epsilon$ satisfies the LDP with good rate function $I:\Xcal \to [0,\infty]$, and let $\phi:\Xcal \to \R$ be any continuous function bounded from above. Then the family $\nu_\epsilon$ defined in \ref{normalized} satisfies the LDP with good rate function
\[J(x)=I(x)-\phi(x)-\inf_{x \in \Xcal}\{I(x)-\phi(x)\}\]
\end{cor}
\begin{proofsect}{Proof}
 Lower Bound: If $\phi:\Xcal \to \R$ is lower semi continuous and $G\subset \Xcal$ open, we define $\phi_G:\Xcal \to \R$ as
\[\phi_G=\begin{cases} \phi(x), \quad x \in G\\
-\infty,\quad \text{otherwise}
\end{cases}\]
then $\phi_G$ is lower semi continuous, because $\{x:\phi_G(x)> a\}=G\cap\{x:\phi(x)> a\}$ which is open for every $a\in \R$. From theorem \ref{vthm} and from lemma \ref{limlo} we have
\[\liminf_{\epsilon \to 0}\epsilon\log \nu_\epsilon(G) =\liminf_{\epsilon \to 0}\epsilon\log\int_G e^{\phi(x)/\epsilon} \mu_\epsilon(\d x)-\liminf_{\epsilon \to 0}\epsilon\log\int_\Xcal e^{\phi(x)/\epsilon} \mu_\epsilon(\d x)\]\[=\liminf_{\epsilon \to 0}\epsilon\log\int_\Xcal e^{\phi_G(x)/\epsilon} \mu_\epsilon(\d x)-\liminf_{\epsilon \to 0}\epsilon\log\int_\Xcal e^{\phi(x)/\epsilon} \mu_\epsilon(\d x)\]
\[\geq\sup_{x\in \Xcal}\{\phi_G(x)-I(x)\}-\sup_{x\in \Xcal}\{\phi(x)-I(x)\}=-\inf_{x \in G}J(x)\]
Upper Bound: If $\phi:\Xcal \to \R$ is upper semi continuous and $F\subset \Xcal$ closed, we define $\phi_F:\Xcal \to \R$ as
\[\phi_F=\begin{cases} \phi(x), \quad x \in F\\
-\infty,\quad \text{otherwise}
\end{cases}\]
then $\phi_F$ is upper semi continuous, because $\{x:\phi_F(x) \geq a\}=F\cap\{x:\phi(x) \geq a\}$ which is closed for every $a\in \R$. From theorem \ref{vthm} and from lemma \ref{lemup} we have
\[\limsup_{\epsilon \to 0}\epsilon\log \nu_\epsilon(F) =\limsup_{\epsilon \to 0}\epsilon\log\int_F e^{\phi(x)/\epsilon} \mu_\epsilon(\d x)-\limsup_{\epsilon \to 0}\epsilon\log\int_\Xcal e^{\phi(x)/\epsilon} \mu_\epsilon(\d x)\]\[=\limsup_{\epsilon \to 0}\epsilon\log\int_\Xcal e^{\phi_F(x)/\epsilon} \mu_\epsilon(\d x)-\limsup_{\epsilon \to 0}\epsilon\log\int_\Xcal e^{\phi(x)/\epsilon} \mu_\epsilon(\d x)\]
\[\leq\sup_{x\in \Xcal}\{\phi_F(x)-I(x)\}-\sup_{x\in \Xcal}\{\phi(x)-I(x)\}=-\inf_{x \in F}J(x)\]
J is a good rate function: Obviously $J$ is non negative. We suppose that there exist a $x_0\in \Xcal$ such that $I(x_0)-\phi(x_0)<\infty$.  For every $y\in \Xcal$ have that $I(y)-\phi(y)\geq -\sup \phi:=-b$. If $c=\inf_{x \in \Xcal}\{I(x)-\phi(x)\}$ then 
\[I(x_0)-\phi(x_0)\geq c \geq -b\]
So if $\phi$ is bounded, then $c\in \R$. Furthermore, since $I$ and $-\phi$ are lower semi continuous, then $J$ is too. Finally, for $a > 0$
\[\{x: J(x)\leq a\}\subset \{x: I(x) \leq a+b+c\}\]
since I is good rate function, the latter set is compact, so J is good rate function as well
\qed
\end{proofsect}

\subsection{Schr\"odinger processes}\label{sec-sc}\bigskip\\

\noindent Here we refer to the paper [FG97] and the book [F88].
\begin{defn}\label{sec-scr}
    Let $(\Omega,\Fcal,(\Fcal_t)_{t\geq 0},\P)$ be a probability space and $S$ be a state space, which we assume that is Polish, and let $(X_t)_{t\ge 0}$ (with $\beta >0$) be a Markov process with underlying probability measure $\P$, and state space $S$. Moreover, let $\mu \in \Mcal_1(S\times S)$ be the joint initial and terminal distribution of $\P$ starting from $t=0$ and terminating at time $t=\beta$. Define the measure 
\[\d \nu (x,y) :=\varphi(x,y) \d \mu (x,y)\] where $\varphi: S\times S \to (0,\infty)$ is a positive, normalized under the measure $\mu $ function. Finally let
\end{defn}
\begin{equation}\label{phisch}
 \d Q:=\varphi (X_0,X_\beta) \d \P
\end{equation}
The process $(X_t)$ under $Q$ is said to be a Schr\"odinger process if $(X_t)$ has the Markov property under the measure $Q$.

\subsubsection{Background.} Let $X_t$ above be the Brownian motion on $S=\R ^d$ with initial distribution $\mu_0$ (that is, $\P$ is the Wiener measure with initial distribution $\mu_0$) and let $P_{s,t}(x,\d y)$, $x\in \R ^d$, $0\leq s<t\leq \beta$ be its transition probabilities. Schr\"odinger initiated the study of processes of the form of $Q$ by considering 
the following problem of Large Deviations. Let $\nu_0$ and $\nu_\beta$ be distributions on $\R ^ d$, and look for the most likely behavior of a large collection of independent Brownian motions governed by $\P$ under the constraint that the empirical distributions at times $t = 0$ and $t = \beta$ are close to $\nu_0$ and to $\nu_\beta$. In the limit of an infinite particle system, this behavior can be described by independent motions governed by $Q$ where the measure $\nu$ minimizes the relative entropy $H(\nu|\mu)$ under the constraints that $\nu_0$ and $\nu_\beta$ are the marginal distributions of $\nu$. See (1.3 page 161 [F88]) for the proof. The density of the entropy minimizing measure admits a factorization of the form 
\begin{equation}\label{factorization}
\frac{\d \nu}{\d \mu}(x,y)=f(x)g(y),
\end{equation}
where $f,g$ are measurable, non negative functions on $\R ^d$. That implies that $Q$ has the structure of $h$-path process, i.e., it is the integral of a function $h$ times the transition probabilities ) where $h$ is a space-time harmonic function of Brownian motion in the sense that $h(s,x)=\int h(t,y)P_{s,t}(x,\d y)$. In particular, we obtain the Markov property and $Q$ is a Schr\"odinger process.

\noindent For a measure $\nu \in \Mcal_1(\R ^d \times \R ^d)$ $\nu \approx \mu$, with $H(\nu|\mu)<\infty$, the Markov process of the associated measure $Q$, is equivalent both to the factorization of $\varphi$ of the form \eqref{factorization} and the minimization of $H(\nu|\mu)$. Such implications are also valid if we replace the Brownian motion on $\R ^d$ with a general Markov process, but not all of them. The difficult part is to derive the factorization \eqref{factorization}. 

\noindent Let $(X_t)$ be a Markov process with undrelying measure $\P$, and $Q$ be a measure with density $\varphi$ as in definition \ref{phisch}. Then for $\t \in [0,\beta]$, we have that 
\[\frac{\d P}{\d Q}\Big|_{\Fcal_t}=\E_{\P}[\phi(X_0,X_\beta)|\Fcal_t]=\int_{S}\phi(X_0,y)p(t,\beta,X_t,\d y):=\varphi_{0,t}(X_0,X_t),\]
where for $x\in S$ we denote by $E^{x}_{\P}$ the expectation with respect to $\P$ when the process starts from $x$, and for $s,t\in [0,\beta]$, with $t>s$, $ x\in S, A\in \Bcal(S)$ we denote by $p(s,t,x,A)$ the transition probabilities of $\P$. 

\noindent Fix $t,s \in [0,\beta]$ with $s<t$. Then we have for every continuous and bounded function $h:S\to \R$ that 
\[\E_Q[h(X_t)|\Fcal_s]=\E_{\P}[h(X_t)\frac{\varphi_{0,t}(X_0,X_t)}{\varphi_{0,s}(X_0,X_s)}|\Fcal_s]\]
To see that, choose a random element of S, $Y$ $\Fcal_s$ measurable. Then
\[\E_Q[h(X_t)Y]=\E_{\P}[Yh(X_t)\varphi_{0,t}(X_0,X_t)]=\]
\[\E_{\P}[Y\E_{\P}[h(X_t)\frac{\varphi_{0,t}(X_0,X_t)}{\varphi_{0,s}(X_0,X_s)}|\Fcal_s]\varphi_{0,s}(X_0,X_s)]=\E_Q[Y\E_{\P}[h(X_t)\frac{\varphi_{0,t}(X_0,X_t)}{\varphi_{0,s}(X_0,X_s)}|\Fcal_s]].\]
So we deduce 
\[Q(X_t\in \d y|\Fcal_s)=\frac{1}{\varphi_{0,s}(X_0,X_s)}\varphi_{0,t}(X_0,y)p(t-s,X_s,\d y).\]
Thus a prediction given $\Fcal_s$ in general will involve both $X_s$ and $X_0$, so $Q$ is not necessarily a Markov process.
However, if $\phi (x,y)=f(x)g(y)$, where f,g are positive, measurable functions on S, then we show that $Q$ is a Schr\"odinger process. Define the function $h:S\times [0,\beta]\to \R$ as
\[h(x,t):=\int_Sg(y)p(1-t,x,\d y),\]
then, $h$ is space-time harmonic in the sense that 
\[h(x,s)=\int_Sh(y,t)p(t-s,x,\d y)\]
for $0\leq s<t\leq \beta$. The density of Q with respect to $\P$ on $\Fcal_t$ can be expressed as
\[\varphi_{0,t}(X_0,X_t)=f(X_0)h(X_t,t)\]
where $\mu_0$ and $\nu_0$ are the initial distributions of $\P$ and Q respectively. This shows that $Q$ has the structure of an h-path process, in particular, Q is a Markov process since h-path transforms preserve the Markov property.


\begin{thebibliography}{WWW98} 


\bibitem[AK08]{AK08}  {\sc S.~Adams,W.~K\"onig},
\newblock Large deviations for symmetrised Brownian bridges
\newblock Prob. Theory Relat. Fields. in press, published online DOI 10.1007/s00440-007-0099-5

\smallskip

\bibitem[B68]{B68}
{\sc P.~Billingsley} ,
\newblock Convergence of probability measures 2nd edition
\newblock {\it Wiley Series in Probability and Statistics}  (1968).
\smallskip

\bibitem[DZ98]{DZ98}
{\sc A.~Dembo} and {\sc O.~Zeitouni},
\newblock {\it Large Deviations Techniques and Applications,}
\newblock 2nd ed., Springer, New York (1998).
\smallskip

\bibitem[FG97]{FG97} 
{\sc H.~F\"ollmer} and {\sc N.~Gantert},
\newblock Entropy minimization and Schr\"odinger processes in infinite dimensions,
\newblock {\it Ann. Probab.} {\bf 25:2}, 901--926 (1997).
\smallskip


\bibitem[CZ95]{CZ95} 
{\sc K.~Chung} and {\sc Z.~Zhao},
\newblock From Brownian motion to Schr\"odinger’s equation,
\newblock Springer Berlin (1995).
\smallskip

\bibitem[MU11]{MU11} 
{\sc L.M~Morato} and {\sc S.~Ugolini},
\newblock Stochastic Description of a Bose–Einstein
Condensate,
\newblock Ann. Henri Poincar´e 12 (2011), 1601–1612.
\smallskip

\bibitem[AD08]{AD08} 
{\sc T.~Dorlas} and {\sc S.~Adams},
\newblock Asymptotic Feynman–Kac formulae for large symmetrised
systems of random walks,
\newblock 
www.imstat.org/aihp
Annales de l’Institut Henri Poincaré - Probabilités et Statistiques
2008, Vol. 44, No. 5, 837–875.
\smallskip


\bibitem[ABK06a]{ABK06a} 
{\sc W.~K\"onig}, {\sc J.B.~Bru} and {\sc S.~Adams},
\newblock Large systems of path-repellent Brownian motions in a
trap at positive temperature,
\newblock 
Electronic journal of Probability Vol. 11 (2006), Paper no. 19, pages 460–485.
\smallskip

\bibitem[RY99]{RY99}
{\sc D.~Revuz} and {\sc M.~Yor},
\newblock {\it Continuous Martingales and Brownian Motion,}
\newblock Springer Berlin (1999).

\smallskip

\bibitem[K17]{K17} {\sc M.~Kohlmann}, \textit{Schr\"odinger Operators and their Spectra}, Lecture notes (2017). 


\smallskip

\bibitem[S98]{S98} {\sc A.~Alain-Sol Sznitman}, \textit{Brownian Motion, Obstacles and Random Media}, Springer Berlin, Heidelberg (1998). 


\smallskip

\bibitem[RS78IV]{RS78IV} {\sc M.~Reed. and B.~Simon. }, \textit{Methods of Modern Mathematical Physics}, Elsevier Science (1978). 


\smallskip

\end{thebibliography}
\end{document}